\documentclass[10pt, reqno]{amsart}
\usepackage{enumerate}
\usepackage{amsmath,amssymb}
\newtheorem{pdef}{Definition}[section] %
\newtheorem{thm}[pdef]{Theorem}        
\newtheorem{cor}[pdef]{Corollary}
\newtheorem{lem}[pdef]{Lemma}
\newtheorem{prop}[pdef]{Proposition}

\newcounter{equationnumber}
\renewcommand{\theequation}{\thesection.\arabic{equation}}
\def\mathletters{
    \addtocounter{equation}{1}
    \edef\@currentlabel{\theequation}
    \setcounter{equationnumber}{\value{equation}}
    \setcounter{equation}{0}
    \edef\theequation{\@currentlabel\noexpand\alph{equation}}
    }

\title{Supports, regularity, and $\boxplus$-infinite \\ divisibility for measures of the
form $(\mu^{\boxplus p})^{\uplus q}$}
\author{Hao-Wei Huang}

\address{Department of Mathematics, Indiana University, 831 East 3rd Street,
Bloomington, IN 47405} \email{huang39@indiana.edu}

\begin{document}

\maketitle\pagestyle{myheadings} \thispagestyle{plain}
\markboth{Hao-Wei Huang}{Supports, regularity, and
$\boxplus$-infinite divisibility for measures of the form
$(\mu^{\boxplus p})^{\uplus q}$}

\begin{abstract} Let $\mathcal{M}$ be the set of Borel probability measures on $\mathbb{R}$.
We denote by $\mu^{\mathrm{ac}}$ the absolutely continuous part of
$\mu\in\mathcal{M}$. The purpose of this paper is to investigate the
supports and regularity for measures of the form $(\mu^{\boxplus
p})^{\uplus q}$, $\mu\in\mathcal{M}$, where $\boxplus$ and $\uplus$
are the operations of free additive and Boolean convolution on
$\mathcal{M}$, respectively, and $p\geq1$, $q>0$. We show that for
any $q$ the supports of $((\mu^{\boxplus p})^{\uplus
q})^{\mathrm{ac}}$ and $(\mu^{\boxplus p})^{\mathrm{ac}}$ contain
the same number of components and this number is a decreasing
function of $p$. Explicit formulas for the densities of
$((\mu^{\boxplus p})^{\uplus q})^{\mathrm{ac}}$ and criteria for
determining the atoms of $(\mu^{\boxplus p})^{\uplus q}$ are given.
Based on the subordination functions of free convolution powers, we
give another point of view to analyze the set of
$\boxplus$-infinitely divisible measures and provide explicit
expressions for their Voiculescu transforms in terms of free and
Boolean convolutions.
\end{abstract} 

\section{Introduction}

For measures $\mu$ and $\nu$ in $\mathcal{M}$, the measure
$\mu\boxplus\nu$ is the free (additive) convolution of $\mu$ and
$\nu$. Thus, $\mu\boxplus\nu$ is the distribution of $X+Y$, where
$X$ and $Y$ are free random variables with distributions $\mu$ and
$\nu$, respectively. Denote by $\phi_\mu$ the Voiculescu transform
of $\mu$ which satisfies the identity
$\phi_{\mu\boxplus\nu}=\phi_\mu+\phi_\nu$ in some truncated cone in
the upper half-plane $\mathbb{C}^+$.

For $n\in\mathbb{N}$, the $n$-fold free convolution
$\mu\boxplus\cdots\boxplus\mu$ is denoted by $\mu^{\boxplus n}$. It
was shown in [\ref{NS}] that the discrete semigroup $\{\mu^{\boxplus
n}:n\in\mathbb{N}\}$ can be embedded in a continuous family
$\{\mu^{\boxplus p}:p\geq1\}$ which satisfies $\mu^{\boxplus
p_1}\boxplus\mu^{\boxplus p_2}=\mu^{\boxplus(p_1+p_2)}$,
$p_1,p_2\geq1$. Any measure in this family satisfies
$\phi_{\mu^{\boxplus p}}=p\phi_\mu$ in some truncated cone in
$\mathbb{C}^+$. We refer the reader to [\ref{BB1},\ref{BB2}, and
\ref{HV3}] for complete developments on the existence of this
continuous family. In the full generalization, Belinschi and
Bercovici used the subordination function to construct the measure
$\mu^{\boxplus p}$, $p>1$, and obtained certain regularity
properties. In [\ref{Huang}], an explicit formula for the density of
$(\mu^{\boxplus p})^{\mathrm{ac}}$ was provided and the relation
between the supports of $\mu$ and $\mu^{\boxplus p}$ was analyzed.
As a consequence, the number $n(p)$ of components in the support of
$\mu^{\boxplus p}$ was shown to be a decreasing function of $p$.

An important class of measures in $\mathcal{M}$ is the set of
$\boxplus$-infinitely divisible measures $\mu$. Recall that $\mu$ is
$\boxplus$-infinitely divisible if for any $n\in\mathbb{N}$ there
exists a measure $\mu_n\in\mathcal{M}$ such that $\mu_n^{\boxplus
n}=\mu$. Another operation of convolution is the Boolean convolution
$\uplus$ introduced by Speicher and Woroudi [\ref{Boolean}]. The
connection among free, Boolean, and classical infinite
divisibilities was thoroughly studied by Bercovici and Pata
[\ref{BP}]. An aspect of this connection between infinite
divisibility with respect to $\boxplus$ and $\uplus$ is the Boolean
Bercovici-Pata bijection $\mathbb{B}$.

Another map $\mathbb{B}_t:\mathcal{M}\to\mathcal{M}$ connecting free
and Boolean convolutions is defined by
\[\mathbb{B}_t(\mu)=\left(\mu^{\boxplus(t+1)}\right)^{\uplus\frac{1}{t+1}},\;\;\;\;\;t\geq0,\;\;\mu\in\mathcal{M}.\]
This map introduced by Belinshi and Nica [\ref{BN1}] satisfies
$\mathbb{B}_t\circ\mathbb{B}_s=\mathbb{B}_{t+s}$, $s,t\geq0$. More
importantly, the map $\mathbb{B}_1$ coincides with the Boolean
Bercovici-Pata bijection $\mathbb{B}$. As a result,
$\mathbb{B}_t(\mu)$ is $\boxplus$-infinitely divisible for any
$\mu\in\mathcal{M}$ and $t\geq1$. This led the authors to associate
to each measure $\mu\in\mathcal{M}$ a nonnegative number
$\mathrm{Ind}(\mu)$, which is called $\boxplus$-divisibility
indicator. For instance, the semicircular and Cauchy distributions
have $\boxplus$-divisibility indicators $1$ and $\infty$,
respectively. It was also shown that $\mu$ is $\boxplus$-infinitely
divisible if and only if $\mathrm{Ind}(\mu)\geq1$. For any measure
$\mu\in\mathcal{M}$ with mean zero and unit variance, denote by
$\Phi(\mu)$ the unique measure in $\mathcal{M}$ such that
$E_\mu=G_{\Phi(\mu)}$. Recall that the free Brownian motion started
at $\nu\in\mathcal{M}$ is the process $\{\nu\boxplus\gamma_t:t>0\}$,
where $\gamma_t$ is the centered semicircular distribution of
variance $t$. The connection between this process and the map
$\mathbb{B}_t$ is via the identity
$E_{\mathbb{B}_t(\mu)}=G_{\Phi(\mu)\boxplus\gamma_t}$. These authors
also studied the regularity of measures in
$\mathbb{B}_t(\mathcal{M})$. In [\ref{BN2}], the same authors
studied the map $\mathbb{B}_t$ on the space $\mathcal{D}_c(k)$ of
distributions of $k$-tuples of self-adjoint elements in a
$C^*$-probability space based on moments and combinatorics. As in
[\ref{BN1}], they showed that $\mathbb{B}_1$ is the multi-variable
Boolean Bercovic-Pata bijection and investigated the relation
between $\mathbb{B}_t$ and free Brownian motion. Later, for measures
$\mu,\nu\in\mathcal{D}_c(k)$, Nica [\ref{Nica}] studied the
so-called subordination distribution of $\mu\boxplus\nu$ with
respect to $\nu$, in which a property related to the present paper
is that $(\mu^{\boxplus p})^{\uplus(p-1)/p}$ is
$\boxplus$-infinitely divisible for any $p>1$. For other further
developments on $\mathbb{B}_t$ and the $\boxplus$-divisibility
indicator of the measure $(\mu^{\boxplus p})^{\uplus q}$, we refer
the reader to [\ref{Japan}].

In the present paper, we mainly use the subordination functions for
the $\boxplus$-convolution powers to study $\boxplus$-infinitely
divisible measures. We show that measures of the form
$(\mu^{\boxplus p})^{\uplus q}$ are $\boxplus$-infinitely divisible
for $\mu\in\mathcal{M}$, $p>1$, and $0<q\leq(p-1)/p$. We also
provide explicit formulas for the Voiculescu transforms of
$\boxplus$-infinitely divisible measures. particularly, the compound
free Poisson distribution with the rate $\lambda$ and jump
distribution $\nu\in\mathcal{M}$ is shown to be of the form
$(\nu_1^{\boxplus(\lambda+1)})^{\uplus\lambda/(\lambda+1)}$,
$\nu_1\in\mathcal{M}$, and its $\boxplus$-divisibility indicator is
calculated as well. In the study of the measures with mean zero and
finite variance $\sigma^2$, we reformulate their
$\boxplus$-divisibility indicators in terms of free Brownian motion:
\[\mathrm{Ind}(\mu)=\sup\left\{t\geq0:E_\mu=\sigma^2G_{\nu_t\boxplus\gamma_{t\sigma^2}}\;\;
\mathrm{for}\;\;\mathrm{some}\;\;\nu_t\in\mathcal{M}\right\}.\] As a
consequence of this reformulation, a measure $\nu\in\mathcal{M}$ can
be written as $\nu_1\boxplus\gamma_t$ for some $\nu_1\in\mathcal{M}$
and $t>0$ if and only if $\mathrm{Ind}(\Phi^{-1}(\nu))>0$. Moreover,
we have $\mathrm{Ind}(\mu)>1$ if and only if
$\phi_\mu=\sigma^2G_{\nu\boxplus\gamma_t}$ for some
$\nu\in\mathcal{M}$ and $t>0$. The work [4] provides solid
foundations for the current research and leads us to investigate the
supports and regularity for the measures $\left(\mu^{\boxplus
p}\right)^{\uplus q}$, $p\geq1$, $q>0$. We prove that the nonatomic
parts of this type of measure are absolutely continuous and the
densities are analytic wherever they are positive. More importantly,
the number of components in the support of $((\mu^{\boxplus
p})^{\uplus q})^{\mathrm{ac}}$ is independent of $q$ and a
decreasing function of $p$. Particularly, $(\mu^{\uplus
q})^{\mathrm{ac}}$ contains the same number of components in the
support for any $q>0$ provided that $\mathrm{Ind}(\mu)>0$.

The paper is organized as follows. Section 2 contains definitions
and basic facts in free probability theory. Section 3 provides
complete descriptions about the connections among free, Boolean
convolutions, and $\boxplus$-infinitely divisible measures. Section
4 investigates the set of $\boxplus$-infinitely divisible measures
with mean zero and finite variance. Section 5 contains results about
the supports and regularity for the measures $\left(\mu^{\boxplus
p}\right)^{\uplus q}$, where $p\geq1$ and $q>0$.

\section{Preliminary}

For any complex number $z$ in $\mathbb{C}$, let $\Re z$ and $\Im z$
be the real and imaginary parts of $z$, respectively. Denote by
$\mathbb{C}^+=\{z\in\mathbb{C}:\Im z>0\}$ the complex upper
half-plane. Consider the set $\mathcal{G}$ defined as
\[\mathcal{G}=\left\{G|G:\mathbb{C}^+\to\mathbb{C}^-\;\;\mathrm{is}\;\;\mathrm{analytic}\;\;\mathrm{and}
\;\;\lim_{y\to\infty}iyG(iy)=1\right\}.\] It is known that a
function $G$ is in $\mathcal{G}$ if and only if there exists some
measure $\mu\in\mathcal{M}$ such that $G$ can be written as
\[G(z)=G_\mu(z):=\int_\mathbb{R}\frac{1}{z-s}\;d\mu(s),\;\;\;\;\;z\in\mathbb{C}^+.\]
The function $G_\mu$ is called the Cauchy transform of $\mu$. The
measure $\mu$ can be recovered from $G_\mu$ as the weak limit of the
measures
\begin{equation} \label{inversion}
d\mu_\epsilon(s)=-\frac{1}{\pi}\Im G_\mu(s+i\epsilon)\;ds
\end{equation} as $\epsilon\to0^+$. This is the Stieltjes inversion
formula. Particularly, if $\Im G$ extends continuously to an open
interval containing some point $x\in\mathbb{R}$ then the density of
the absolutely continuous part of $\mu$ at $x$ is given by $-\Im
G(x)/\pi$.

Another class of functions which is closely related to $\mathcal{G}$
and plays a significant role in free probability theory is the
following set
\[\mathcal{F}=\left\{F|F:\mathbb{C}^+\to\mathbb{C}^+\;\;\mathrm{is}\;\;\mathrm{analytic}\;\;\mathrm{and}
\;\;\lim_{y\to\infty}\frac{F(iy)}{iy}=1\right\}.\] A function $F$
belongs to $\mathcal{F}$ if and only if $F=F_\mu:=1/G_\mu$ for some
$\mu\in\mathcal{M}$. The function $F_\mu$ is called the reciprocal
Cauchy transform of $\mu$. Any function $F\in\mathcal{F}$ has the
property $\Im F(z)\geq\Im z$ for $z\in\mathbb{C}^+$ and has a
Nevanlinna representation of the form
\begin{equation}\label{NeF}
F(z)=\Re F(i)+z+\int_\mathbb{R}\frac{1+sz}{s-z}\;d\rho(s),
\end{equation}
where $\rho$ is some finite positive Borel
measure on $\mathbb{R}$. Moreover, the function $F$ has a right
inverse $F_\mu^{-1}$ with respect to composition, which is defined
on the truncated cone
\[\Gamma_{\alpha,\beta}=\{x+iy\in\mathbb{C}^+:|x|\leq\alpha
y,\;|y|\geq\beta\}\] of the upper half-plane for some
$\alpha,\beta>0$. The function
$\phi_\mu:\Gamma_{\alpha,\beta}\to\mathbb{C}^-\cup\mathbb{R}$
defined by
\[\phi_\mu(z)=F_\mu^{-1}(z)-z,\;\;\;\;\;z\in\Gamma_{\alpha,\beta},\]
is called the Voiculescu transform of $\mu$. As indicated in the
introduction, for $\mu,\nu\in\mathcal{M}$ and $z$ in some truncated
cone in $\mathbb{C}^+$ the following identity holds:
\[\phi_{\mu\boxplus\nu}(z)=\phi_\mu(z)+\phi_\nu(z).\] Particularly, the identity
$F_{\mu\boxplus\delta_a}(z)=F_\mu(z-a)$ holds for $z\in\mathbb{C}^+$
and $a\in\mathbb{R}$.

The reciprocal Cauchy transform $F_\mu$ can be used to locate the
atoms of $\mu$. A point $\alpha$ is an atom of $\mu$ if and only if
$F_\mu(\alpha)=0$ (that is, $F_\mu$ is defined and takes the value
$0$ at the point $\alpha$) and the Julia-Carath\'{e}odory derivative
$F_\mu'(\alpha)$ (which is the limit of
\[\frac{F_\mu(z)-F_\mu(\alpha)}{z-\alpha}\] as $z\to\alpha$ nontangentially, i.e., $(\Re
z-\alpha)/\Im z$ stays bounded and $z\in\mathbb{C}^+$) is finite, in
which case $\mu(\{\alpha\})=1/F_\mu'(\alpha)$.

Given any measure $\mu\in\mathcal{M}$, the function
$E_\mu(z)=z-F_\mu(z)$ is called the energy function associated with
$\mu$ and belongs to the following set
\[\mathcal{E}=\left\{E|E:\mathbb{C}^+\to\mathbb{C}^-\cup\mathbb{R}\;\;\mathrm{is}\;\;\mathrm{analytic}\;\;\mathrm{and}
\;\;\lim_{y\to\infty}\frac{E(iy)}{iy}=0\right\}.\] Conversely, any
function $E$ in $\mathcal{E}$ is the energy function of some
$\mu\in\mathcal{M}$ whose Nevanlinna representation is given by
\begin{equation} \label{NeE}
E(z)=\Re E(i)+\int_\mathbb{R}\frac{1+sz}{z-s}\;d\rho(s),
\end{equation} where $\rho$ is some finite positive Borel measure on
$\mathbb{R}$. Observe that we have the inclusion
$\mathcal{G}\subset\mathcal{E}$. Indeed, for any measure
$\mu\in\mathcal{M}$ it was proved in [\ref{Maa}] that $\mu$ has mean
zero and finite variance $\sigma^2$, i.e.,
\[\int_\mathbb{R}s\;d\mu(s)=0\;\;\;\;\;\mathrm{and}\;\;\;\;\;\int_\mathbb{R}s^2\;d\mu(s)=\sigma^2\]
if and only if there exists some unique $\nu\in\mathcal{M}$ such
that
\begin{equation} \label{EMaa} E_\mu=\sigma^2G_\nu.
\end{equation} If $\sigma^2=1$, let $\Phi(\mu)$ be the unique measure
satisfying $E_\mu=G_{\Phi(\mu)}$. The Eq. (\ref{EMaa}) particularly
shows that $\mu^{\uplus1/\sigma^2}$ has mean zero and unit variance,
i.e., $E_\mu=\sigma^2G_{\Phi(\mu^{\uplus1/\sigma^2})}$.

Next, consider the set
\[\mathcal{H}=\left\{H|H:\mathbb{C}^+\to\mathbb{C}\;\;\mathrm{is}\;\;\mathrm{analytic},\;\Im H(z)\leq\Im
z, \;\;z\in\mathbb{C}^+,\;\;\mathrm{and}
\;\;\lim_{y\to\infty}\frac{H(iy)}{iy}=1\right\},\] which plays an
important role in the investigation of the free convolution powers
of measures in $\mathcal{M}$. Indeed, for any $H\in\mathcal{H}$ the
function $2z-H(z)\in\mathcal{F}$ is the reciprocal Cauchy transform
of some measure in $\mathcal{M}$. More importantly, the right
inverses of the functions in $\mathcal{H}$ can be used to construct
the $p$-th $\boxplus$-convolution power $\mu^{\boxplus p}$,
$p\geq1$, of any measure $\mu\in\mathcal{M}$. We list below the
properties needed in this paper. For more details, we refer the
reader to [\ref{BB1},\ref{BB2}, and \ref{Huang}].

\begin{prop} \label{prop2.1} For any $\mu\in\mathcal{M}$ and $p>1$,
define the function
\[H_p(z)=pz+(1-p)F_\mu(z),\;\;\;\;\;z\in\mathbb{C}^+,\] the set
$\Omega_p=\{z\in\mathbb{C}^+:\Im H_p(z)>0\}$, and the function
$f_\mu:\mathbb{R}\to\mathbb{R}_+\cup\{\infty\}$ as
\begin{equation} \label{fmu}
f_\mu(x)=\int_\mathbb{R}\frac{s^2+1}{(s-x)^2}\;d\rho(s),\;\;\;\;\;x\in\mathbb{R},
\end{equation} where $\rho$ is the measure in the Nevanlinna representation $(\ref{NeF})$
of $F_\mu$.
\begin{enumerate} [$\qquad(1)$]
\item {The function $H_p$ is in $\mathcal{H}$ and the set $\Omega_p$ is
a simply connected domain whose boundary is the graph of the
continuous function $f_p:\mathbb{R}\to[0,\infty)$, where
\[f_p(x)=\inf\left\{y>0:\frac{\Im
E_\mu(x+iy)}{y}>\frac{-1}{p-1}\right\},\;\;\;\;\;x\in\mathbb{R}.\]}
\item {For $x\in\mathbb{R}$, $f_p(x)=0$ if and only if
$f_\mu(x)\leq1/(p-1)$, while $z\in\Omega_p$ if and only if
\[\int_\mathbb{R}\frac{s^2+1}{|s-z|^2}\;d\rho(s)<\frac{1}{p-1}.\]
Consequently, the functions $E_\mu$ and $H_p$ have continuous
extensions to $\overline{\Omega_p}$ which are Lipschitz continuous
with Lipschitz constants $1/(p-1)$ and $2$, respectively. Moreover,
Eq. $(\ref{NeF})$ holds for $z\in\overline{\Omega_p}$.}
\item {There exists an analytic function
$\omega_p:\mathbb{C}^+\to\mathbb{C}^+$ extending continuously to
$\mathbb{C}^+\cup\mathbb{R}$ such that $H_p(\omega_p(z))=z$ holds
for $z\in\mathbb{C}^+\cup\mathbb{R}$. Consequently,
$\Omega_p=\omega_p(\mathbb{C}^+)$, $\omega_p(H_p(z))=z$ holds for
$z\in\overline{\Omega_p}$, and
\[\frac{|z_1-z_2|}{2}\leq|\omega_p(z_1)-\omega_p(z_2)|,\;\;\;\;\;z_1,z_2\in\mathbb{C}^+\cup\mathbb{R}.\]}
\item {The function $\omega_p$ is analytic in a neighborhood of $x$
wherever $\omega_p(x)\not\in\mathbb{R}$.}
\item {Let
$\mu^{\boxplus p}$ be the unique measure in $\mathcal{M}$ whose
reciprocal Cauchy transform satisfies
\begin{equation} \label{p-power}
F_{\mu^{\boxplus
p}}(z)=\frac{p\omega_p(z)-z}{p-1},\;\;\;\;\;z\in\mathbb{C}^+.
\end{equation} Then there exist some $\alpha,\beta>0$ such that
\[\phi_{\mu^{\boxplus
p}}(z)=p\phi_\mu(z),\;\;\;\;\;z\in\Gamma_{\alpha,\beta}.\] Moreover,
the function $\omega_p$ is the subordination function of
$\mu^{\boxplus p}$ with respect to $\mu$, i.e.,
\[F_{\mu^{\boxplus
p}}(z)=F_\mu(\omega_p(z)),\;\;\;\;\;z\in\mathbb{C}^+\cup\mathbb{R},\]
and consequently
\[F_{\mu^{\boxplus
p}}(H_p(z))=F_\mu(z),\;\;\;\;\;z\in\overline{\Omega_p}.\]}
\end{enumerate}
\end{prop}

Complete characterizations of the supports of $\mu^{\boxplus p}$
were given in [\ref{Huang}]. Following the notations in Proposition
\ref{prop2.1}, we give below the results needed in the current
research.

\begin{thm} \label{Hthm} For $\mu\in\mathcal{M}$, define the function
$\psi_p:\mathbb{R}\to\mathbb{R}$ by $\psi_p(x)=H_p(x+if_p(x))$,
$x\in\mathbb{R}$, and the set $V_p^+=\{x\in\mathbb{R}:f_p(x)>0\}$.
Then the following statements are true.
\begin{enumerate} [$\qquad(1)$]
\item {The function $\psi_p$ is a homeomorphism on $\mathbb{R}$.}
\item {The measure $(\mu^{\boxplus
p})^{\mathrm{ac}}$ is concentrated on the set
$\psi_p(\overline{V_p^+})$ with density
\begin{equation} \label{density}
\frac{d(\mu^{\boxplus
p})^{\mathrm{ac}}}{dx}(\psi_p(x))=\frac{(p-1)pf_p(x)}{\pi|px-\psi_p(x)+ipf_p(x)|},
\;\;\;\;\;x\in V_p^+.
\end{equation} }
\item {The number of the components in the support of $\mu^{\boxplus
p}$ is a decreasing function of $p$.}
\end{enumerate}
\end{thm}

The set of $\boxplus$-infinitely divisible measures in $\mathcal{M}$
is closed under weak convergence of probability measures. As shown
in [\ref{HV2}], a necessary and sufficient condition for $\mu$ to be
$\boxplus$-infinitely divisible is that $\phi_\mu$ belong to
$\mathcal{E}$.

The Boolean convolution introduced in [\ref{Boolean}] was defined
via the functions in $\mathcal{E}$. Given $\mu_1$ and $\mu_2$ in
$\mathcal{M}$, the measure $\nu$ satisfying the relation
\[E_\nu=E_{\mu_1}+E_{\mu_2}\] is called the Boolean convolution of
$\mu_1$ and $\mu_2$, and it is denoted $\mu_1\uplus\mu_2$. For
$\mu\in\mathcal{M}$ and a positive integer $n$, the $n$-fold Boolean
convolution $\mu\uplus\cdots\uplus\mu$ denoted by $\mu^{\uplus n}$
satisfies $E_{\mu^{\uplus n}}=nE_\mu$. This can be extended
naturally to the case when the exponent $n$ is not an integer. That
is, for every $q\geq0$ the $q$-th $\uplus$-convolution power
$\mu^{\uplus q}$ is defined as the unique measure in $\mathcal{M}$
satisfying
\[E_{\mu^{\uplus q}}=qE_\mu.\]

The following theorem builds the connection between
$\boxplus$-infinitely divisible measures and the Boolean
convolution, which was thoroughly investigated in [\ref{BP}].

\begin{thm} \label{thm2.2} Let $\{\mu_n\}$ be a sequence in
$\mathcal{M}$ and let $k_1<k_2<\cdots$ be a sequence of positive
integers. Then the following statements $(1)$-$(3)$ are equivalent:
\begin{enumerate} [$\qquad(1)$]
\item {$\mu_n^{\boxplus k_n}\to\mu_{\boxplus}^{c,\rho}$ weakly as
$n\to\infty$;} \item {$\mu_n^{\uplus k_n}\to\mu_{\uplus}^{c,\rho}$
weakly as $n\to\infty$;} \item {the measures
\[k_n\frac{s^2}{s^2+1}\;d\mu_n(s)\to d\rho(s)\] weakly as $n\to\infty$ and
\[\lim_{n\to\infty}\int_\mathbb{R}\frac{k_ns}{s^2+1}\;d\mu_n(s)=c.\]}
\end{enumerate}
If $(1)$-$(3)$ hold then $\mu_{\boxplus}^{c,\rho}$ is
$\boxplus$-infinitely divisible and
\[\phi_{\mu_{\boxplus}^{c,\rho}}(z)=E_{\mu_{\uplus}^{c,\rho}}(z)=c+\int_\mathbb{R}\frac{1+sz}{z-s}\;d\rho(s),
\;\;\;\;\;z\in\mathbb{C}^+.\]
\end{thm}

For $\mu\in\mathcal{M}$, Theorem \ref{thm2.2} shows that
$(\mu^{\uplus1/n})^{\boxplus n}$ converges weakly to some
$\boxplus$-infinitely divisible measure $\mathbb{B}(\mu)$ satisfying
$\phi_{\mathbb{B}(\mu)}=E_\mu$. Conversely, for any
$\boxplus$-infinitely divisible measure $\nu$ the sequence
$(\nu^{\boxplus1/n})^{\uplus n}$ converges weakly to some $\mu$
satisfying $\phi_\nu=E_\mu$. Since $E$ determines the measure
uniquely, the map $\mathbb{B}$ induces a bijective map from
$\mathcal{M}$ onto the set of $\boxplus$-infinitely divisible
measures. This map $\mathbb{B}$ is called the Boolean Bercovici-Pata
bijection, which coincides with $\mathbb{B}_1$ as indicated in the
introduction.

In the study of $\boxplus$-infinitely divisible measures, there is
one useful tool introduced in [\ref{BN1}] called
$\boxplus$-divisibility indicator:
\[\mathrm{Ind}(\mu)=\sup\{t\geq0:\mu\in\mathbb{B}_t(\mathcal{M})\},\;\;\;\;\;\mu\in\mathcal{M}.\]
Any measure $\mu\in\mathcal{M}$ with finite support has
$\mathrm{Ind}(\mu)=0$, while $\mu$ is $\boxplus$-infinitely
divisible if and only if $\mathrm{Ind}(\mu)\geq1$. In general, for
any $t\geq0$ and $\mu\in\mathcal{M}$ we have
\begin{equation} \label{Nindicator}
\mathrm{Ind}(\mathbb{B}_t(\mu))=t+\mathrm{Ind}(\mu)
\end{equation} For
$\boxplus$-divisibility indicators of some specific measures, we
refer the reader to [\ref{BN1}].

\setcounter{equation}{0}
\section{$\boxplus$-infinite divisibility of $\mathbb{B}_{p,q}(\mu)$}

For any $\mu\in\mathcal{M}$ and $p\geq1$, $q>0$, denote
$\mathbb{B}_{p,q}(\mu)=\left(\mu^{\boxplus p}\right)^{\uplus q}$.
Particularly, $\mathbb{B}_{t+1,1/(t+1)}=\mathbb{B}_t$ for any
$t\geq0$ and $\mathbb{B}_{2,1/2}=\mathbb{B}$. In this section, we
mainly use Proposition \ref{prop2.1} to investigate the measure
$\mathbb{B}_{p,q}(\mu)$. Throughout the paper, the number $r^*$
stands for the conjugate exponent of any number $r>0$, i.e.,
\[\frac{1}{r}+\frac{1}{r^*}=1,\;\;\;\;\;r\neq1\] and $r^*=\infty$ if $r=1$.
Note that we have $r^*<0$ if $r\in(0,1)$.

For $p>1$, by (\ref{p-power}) and the definition of Boolean
convolution power we have
\begin{equation} \label{general}
F_{\mathbb{B}_{p,q}(\mu)}(z)=\frac{pq\omega_p(z)-(1+pq-p)z}{p-1},\;\;\;\;\;z\in\mathbb{C}^+\cup\mathbb{R}.
\end{equation}
As a special case of (\ref{general}), if $1+pq-p=0$, i.e., $q=1/p^*$
then
\[F_{\mathbb{B}_{p,1/p^*}(\mu)}(z)
=\frac{p\omega_p(z)}{p^*(p-1)}=\omega_p(z),\;\;\;\;\;z\in\mathbb{C}^+\cup\mathbb{R}.\]
This yields that the Voiculescu transform
$\phi_{\mathbb{B}_{p,1/p^*}(\mu)}$ of the measure
$\mathbb{B}_{p,1/p^*}(\mu)$ has an analytic continuation to
$\mathbb{C}^+$, which is given by
\begin{align*}
\phi_{\mathbb{B}_{p,1/p^*}(\mu)}(z)&=F_{\mathbb{B}_{p,1/p^*}(\mu)}^{-1}(z)-z=H_p(z)-z \\
&=E_{\mu^{\uplus(p-1)}}(z),\;\;\;\;\;z\in\mathbb{C}^+.
\end{align*}
These observations are recorded in the following result.

\begin{prop} \label{3.1} For any measure $\mu\in\mathcal{M}$ and number $p>1$,
the measure $\mathbb{B}_{p,1/p^*}(\mu)$ is $\boxplus$-infinitely
divisible, the function $F_{\mathbb{B}_{p,1/p^*}(\mu)}$ extends
continuously to $\mathbb{C}^+\cup\mathbb{R}$,
\[F_{\mu^{\boxplus p}}(z)=F_\mu\left(F_{\mathbb{B}_{p,1/p^*}(\mu)}(z)\right),\;\;\;\;\;z\in\mathbb{C}^+\cup\mathbb{R},\]
and the Voiculescu transform of $\mathbb{B}_{p,1/p^*}(\mu)$ can be
expressed as
\[\phi_{\mathbb{B}_{p,1/p^*}(\mu)}=E_{\mu^{\uplus(p-1)}}.\] In
particular, the above statements hold for $\mathbb{B}_1$.
\end{prop}

Observe that Proposition \ref{3.1} provides an easy way to prove
that $\mathbb{B}_1(\mu)$, $\mu\in\mathcal{M}$, is identically equal
to the image $\mathbb{B}(\mu)$ of $\mu$ under the Boolean
Bercovici-Pata bijection. Indeed, by Theorem \ref{thm2.2} and
Proposition \ref{3.1} we obtain
\begin{equation} \label{bijection}
\phi_{\mathbb{B}(\mu)}=E_\mu=\phi_{\mathbb{B}_1(\mu)}.
\end{equation}
In [\ref{Nica}], results similar to Proposition \ref{3.1} for the
joint distributions for $k$-tuples of selfadjoint elements in a
$C^*$-probability space were obtained by combinatorial tools. We
refer the reader to the same paper for the so-called $k$-tuple
Boolean Bercovici-Pata bijection and related results.

The following lemma contains some basic properties of the map
$\mathbb{B}_{p,q}$ which is frequently used in the sequel. The
identity in (\ref{formula}) can be obtained by [\ref{BN1},
Proposition 3.1]. Here we provide an alternative proof using
Proposition \ref{3.1}.

\begin{lem} If $\mu\in\mathcal{M}$, $p,p_1\geq1$, and $q,q_1>0$
then
\begin{equation} \label{formula}
\left(\mu^{\uplus q}\right)^{\boxplus
p}=\mathbb{B}_{1+pq-q,\frac{pq}{1+pq-q}}(\mu),
\end{equation}
\begin{equation} \label{two}
\mathbb{B}_{p_1,q_1}\circ\mathbb{B}_{p,q}=
\mathbb{B}_{p(1+p_1q-q),\frac{p_1q_1q}{1+p_1q-q}},
\end{equation}
and
\begin{equation} \label{t}
\mathbb{B}_{p,q}=\mathbb{B}_t\circ\mathbb{B}_{p(1-qt),\frac{q}{1-qt}},
\;\;\;\;\;0\leq t\leq\frac{1}{p^*q}.
\end{equation}
\end{lem}

\begin{pf} It suffices to show the lemma for $p>1$.
By Proposition \ref{prop2.1}, we have
\[F_{(\mu^{\uplus q})^{\boxplus p}}(z)=\frac{p\omega(z)-z}{p-1},\;\;\;\;\;z\in\mathbb{C}^+,\]
where the function $\omega$ is the right inverse of the function
\begin{align*}
H(z)&=pz+(1-p)F_{\mu^{\uplus q}} \\
&=(1+pq-q)z+(q-pq)F_\mu(z),\;\;\;\;\;z\in\mathbb{C}^+.
\end{align*}
On the other hand, since the number $1+pq-q>1$ whose conjugate
exponent is
\[(1+pq-q)^*=\frac{1+pq-q}{pq-q},\] by Proposition \ref{3.1} we
see that $\omega=F_\nu$, where
\begin{equation} \label{nu} \nu=\mathbb{B}_{1+pq-q,\frac{pq-q}{1+pq-q}}(\mu).
\end{equation} Then using the definition of the Boolean convolution power and
(\ref{nu}) gives
\[\frac{p\omega(z)-z}{p-1}=\frac{pF_\nu(z)-z}{p-1}=F_{\nu^{\uplus p^*}}(z),\;\;\;\;\;z\in\mathbb{C}^+,\]
and $\nu^{\uplus p^*}=\mathbb{B}_{1+pq-q,\frac{pq}{1+pq-q}}(\mu)$,
whence the formula in (\ref{formula}) follows. The equality in
(\ref{two}) follows directly from (\ref{formula}). Finally, note
that if $t\in[0,1/(p^*q)]$ then
\[p(1-qt)\geq1\;\;\;\;\;\mathrm{and}\;\;\;\;\;\frac{q}{1-qt}>0,\] whence the measure
$\mathbb{B}_{p(1-qt),q/(1-qt)}(\mu)$ is defined and (\ref{t}) holds
by (\ref{two}).
\end{pf} \qed

If $p>1$ and $0<q<1/p^*$ (or, equivalently, $1+pq-p<0$) then
(\ref{t}) yields
$\mathbb{B}_{p,q}=\mathbb{B}_1\circ\mathbb{B}_{p_1,q_1}$, where
\[p_1=p(1-q)>1\;\;\;\;\;\mathrm{and}\;\;\;\;\;q_1=\frac{q}{1-q}>0.\]
Using Proposition \ref{3.1} and the preceding discussions gives the
following result.

\begin{thm} \label{3.3} If $\mu\in\mathcal{M}$, $p>1$, and $0<q\leq1/p^*$ then
the following statements hold.
\begin{enumerate} [$\qquad(1)$]
\item {The measure $\mathbb{B}_{p,q}(\mu)$ is $\boxplus$-infinitely
divisible.}
\item {For any $n\in\mathbb{N}$,
\[{\mathbb{B}_{p,q}(\mu)}=
\left({\mathbb{B}_{\frac{p(n+q-nq)}{n},\frac{q}{n(1-q)+q}}(\mu)}\right)^{\boxplus
n}.\]}
\item {The Voiculescu transform of $\mathbb{B}_{p,q}(\mu)$ can be
expressed as
\[\phi_{\mathbb{B}_{p,q}(\mu)}=E_{\mathbb{B}_{p_1,q_1}(\mu)}.\]}
\item {For $r>0$, let $\nu_r=\mathbb{B}_{p_1,r}(\mu)$. Then
\[F_{\mathbb{B}_{p,(1-q)r+q}(\mu)}(z)
=F_{\nu_r}\left(F_{\mathbb{B}_{p,q}(\mu)}(z)\right),\;\;\;\;\;z\in\mathbb{C}^+\cup\mathbb{R}.\]}
\end{enumerate}
Particularly, for any $t\geq1$ the measure $\mathbb{B}_t(\mu)$ is
$\boxplus$-infinitely divisible,
\[\phi_{\mathbb{B}_t(\mu)}=E_{\mathbb{B}_{t-1}(\mu)},\] and
$F_{\mathbb{B}_t(\mu)}$ is the subordination function of
$\mu^{\boxplus(t+1)}$ with respect to $\mu^{\boxplus t}$, that is,
\[F_{\mu^{\boxplus(t+1)}}(z)=F_{\mu^{\boxplus t}}(F_{\mathbb{B}_t(\mu)}(z)),
\;\;\;\;\;z\in\mathbb{C}^+\cup\mathbb{R}.\]
\end{thm}

\begin{pf} The assertions (1) and (3) were proved (particularly, $p=t+1$ and $q=(t+1)^{-1}$
satisfy the condition $1+pq-p\leq0$ if and only if $t\geq1$). Next,
observe that
\[\frac{p(n+q-nq)}{n}>p-pq\geq1\;\;\;\;\;\mathrm{and}\;\;\;\;\;
\frac{q}{n(1-q)+q}>0,\] whence the assertion (2) follows from
(\ref{two}). By (3), $F_{\mathbb{B}_{p,q}(\mu)}^{-1}$ can be
expressed as
\begin{equation} \label{1}
F_{\mathbb{B}_{p,q}(\mu)}^{-1}(z)=
\left(1+\frac{q_1}{r}\right)z-\frac{q_1}{r}F_{\nu_r}(z),\;\;\;\;\;z\in\mathbb{C}^+.
\end{equation} Since
$\nu_r^{\boxplus\left(1+q_1/r\right)}=\mathbb{B}_{p,(1-q)r+q}(\mu)$
by (\ref{formula}), Proposition \ref{prop2.1}(4) and (\ref{1}) imply
the assertion (4). Letting $r=1$ in (4) yields the last assertion.
\end{pf} \qed

Observe that if $\mu$ is $\boxplus$-infinitely divisible then
$\mu\in\mathbb{B}(\mathcal{M})$, i.e., the measure
$\mathbb{B}^{-1}(\mu)=(\mu^{\uplus2})^{\boxplus1/2}$ is defined. In
order to investigate the measure of the form $\mu^{\boxplus p}$,
$0<p<1$ (that is, $\mu=\nu^{\boxplus1/p}$ for some
$\nu\in\mathcal{M}$ ), we need the following lemma. This lemma was
also provided in [\ref{Japan}]; however, the case
$\mathrm{Ind}(\mu)=\infty$ was not considered there.

\begin{lem} \label{3.4} For any measure $\mu\in\mathcal{M}$ and any number $q>0$, we have
\[\mathrm{Ind}\left(\mu^{\uplus q}\right)=\frac{\mathrm{Ind}(\mu)}{q}.\]
\end{lem}

\begin{pf} First claim the inequality $\mathrm{Ind}(\mu^{\uplus q})\geq\mathrm{Ind}(\mu)/q$ holds.
It clearly holds if $\mathrm{Ind}(\mu)=0$. Next, consider the case
$\mathrm{Ind}(\mu)>0$. Then for any finite $r$ with
$0<r<\mathrm{Ind}(\mu)$ pick a measure $\nu\in\mathcal{M}$ such that
$\mu=\mathbb{B}_r(\nu)$, from which we obtain that
\[\mathrm{Ind}\left(\mu^{\uplus q}\right)=\mathrm{Ind}\left(\left(\nu^{\boxplus(r+1)}\right)^{\uplus
\frac{q}{r+1}}\right)=\mathrm{Ind}\left(\mathbb{B}_{r/q}(\nu_1)\right)\geq
\frac{r}{q},\] where $\nu_1\in\mathcal{M}$ and (\ref{t}) is used in
the second equality above. If $\mathrm{Ind}(\mu)=\infty$ then
letting $r\uparrow\infty$ gives $\mathrm{Ind}(\mu^{\uplus
q})=\infty$ as well, which implies the desired inequality.
Otherwise, letting $r\uparrow\mathrm{Ind}(\mu)<\infty$ yields the
claim. By considering the identity $\mu=\left(\mu^{\uplus
q}\right)^{\uplus1/q}$ and using the claim, we obtain the opposite
inequality, and then the proof is complete.
\end{pf} \qed

Now we are able to determine for what value of $p\in(0,1)$ the
$p$-th $\boxplus$-convolution power of a measure is defined. The
following implication that (2) implies (1) was proved in
[\ref{BN1}].

\begin{prop} \label{3.5} Suppose that $\mu\in\mathcal{M}$ and fix a number $p\in(0,1)$. Then the
following statements $(1)$-$(3)$ are equivalent:
\begin{enumerate} [$\qquad(1)$]
\item {$\mu^{\boxplus p}$ is defined, that is, $\mu=\nu^{\boxplus1/p}$
for some $\nu\in\mathcal{M}$;}
\item {$1-p\leq\mathrm{Ind}(\mu)$;}
\item {$\mu^{\uplus(1-p)}$ is $\boxplus$-infinitely divisible.}
\end{enumerate}
If $(1)$-$(3)$ hold and
$\Omega_p=F_{\mu^{\uplus(1-p)}}(\mathbb{C}^+)$ then the following
statements are true.
\begin{enumerate} [$\qquad(a)$]
\item {The reciprocal Cauchy transform $F_\mu$ extends continuously to
$\mathbb{C}^+\cup\mathbb{R}$.}
\item {The identity $F_{\mu^{\boxplus
p}}\left(F_{\mu^{\uplus(1-p)}}(z)\right)=F_\mu(z)$ holds for
$z\in\mathbb{C}^+\cup\mathbb{R}$.}
\item {The identity $F_{\mu^{\boxplus
p}}(z)=F_\mu\left(F_{\mu^{\uplus(1-p)}}^{-1}(z)\right)$ holds for
$z\in\overline{\Omega_p}$.}
\item {The Voiculescu transform of $\mu^{\uplus(1-p)}$ can be expressed as
\[\phi_{\mu^{\uplus(1-p)}}=E_{\left(\mu^{\boxplus
p}\right)^{\uplus(1/p-1)}}=E_{\mathbb{B}^{-1}(\mu^{\uplus(1-p)})}.\]}
\end{enumerate}
\end{prop}

\begin{pf} If $\nu^{\boxplus1/p}=\mu$ then by Lemma \ref{3.4} we have
\[\mathrm{Ind}(\mu)=p
\mathrm{Ind}\left(\mathbb{B}_{1/p-1}(\nu)\right)\geq
p\left(\frac{1}{p}-1\right)=1-p,\] which shows that (1) implies (2).
Applying Proposition \ref{3.1} to $\nu$ and $1/p$ yields the
assertions (3). If (3) holds, i.e.,
$\phi_{\mu^{\uplus(1-p)}}\in\mathcal{E}$ then the following function
\[F(z)=\frac{-p\phi_{\mu^{\uplus(1-p)}}(z)+(1-p)z}{1-p},\;\;\;\;z\in\mathbb{C}^+,\]
belongs to $\mathcal{F}$, i.e., $F=F_\nu$ for some
$\nu\in\mathcal{M}$, by which $F_{\mu^{\uplus(1-p)}}^{-1}$ can be
written as
\begin{equation} \label{p<1o}
F_{\mu^{\uplus(1-p)}}^{-1}(z)=\frac{1}{p}z+\left(1-\frac{1}{p}\right)F_\nu(z),\;\;\;\;\;z\in\mathbb{C}^+.
\end{equation}
Then by Proposition 2.1 and the definition of the
$\uplus$-convolution power we have
\[F_{\nu^{\boxplus1/p}}(z)=\frac{\frac{1}{p}F_{\mu^{\uplus(1-p)}}(z)-z}{\frac{1}{p}-1}
=F_\mu(z),\;\;\;\;\;z\in\mathbb{C}^+,\] whence (1) holds. The
assertion (a) holds since $\mu=\nu^{\boxplus1/p}$ and $1/p>1$, while
(b)-(d) follow from the preceding discussions, (\ref{bijection}),
and Proposition \ref{prop2.1}(5).
\end{pf} \qed

The proof of the preceding proposition also gives the construction
of the measure $\mu^{\boxplus p}$ whenever it is defined for
$p\in(0,1)$. Indeed, by (\ref{p<1o}) the right inverse $\omega_p$ of
the function $H_p(z)=pz+(1-p)F_\mu(z)$ ($H_p=F_{\mu^{\uplus(1-p)}}$)
satisfies the relation
\begin{equation} \label{p<1} F_{\mu^{\boxplus
p}}(z)=\frac{p\omega_p(z)-z}{p-1},\;\;\;\;z\in\mathbb{C}^+,
\end{equation} and we have $F_{\mu^{\boxplus p}}(z)=F_\mu(\omega_p(z))$ for
$z\in H_p(\mathbb{C}^+)$.

The following proposition can be proved by [\ref{BN1}, Proposition
3.1]. It can be also obtained by using (\ref{formula}),
(\ref{general}), and (\ref{p<1}), and we leave the proof for the
reader.

\begin{prop} \label{3.6} Let $\mu\in\mathcal{M}$ and for $p,q>0$ let $q'=1+pq-p$ and $p'=pq/q'$.
\begin{enumerate} [$\qquad(1)$]
\item {If $\mu^{\boxplus p}$ is defined and $q'>0$ then we have the
following identity
\begin{equation} \label{formula2}\left(\mu^{\boxplus p}\right)^{\uplus q}=\left(\mu^{\uplus q'}\right)^{\boxplus
p'}.
\end{equation}}
\item {The formula $(\ref{formula})$ holds for either
\begin{enumerate} [$(a)$]
\item {$p\geq1$ or}
\item {$1-\mathrm{Ind}(\mu^{\uplus q})\leq p<1$ and $1+pq-q>0$.}
\end{enumerate}}
\end{enumerate}
\end{prop}

It was proved in [\ref{BN1}] that $\mu\in\mathbb{B}_t(\mathcal{M})$
for any finite $t$ with $0\leq t\leq\mathrm{Ind}(\mu)$. In the
following proposition we give an explicit expression for the measure
$\mu_t$ so that $\mu=\mathbb{B}_t(\mu_t)$. The reader should be
aware of that this conclusion holds under the essential condition
that $t$ has to be finite and this condition may not be noticed
without caution.

\begin{prop} If $\mu\in\mathcal{M}$ then for any finite number $t$
with $-1<t\leq\mathrm{Ind}(\mu)$ there exists a  unique measure
$\mu_t\in\mathcal{M}$ such that
\begin{equation} \label{N}
\mu=\left(\mu_t^{\boxplus(1+t)}\right)^{\uplus\frac{1}{1+t}},
\end{equation}
in which case $\mu_t$ can be expressed as
\begin{equation} \label{N1}
\mu_t=\left(\mu^{\uplus(t+1)}\right)^{\boxplus\frac{1}{t+1}}.
\end{equation} Particularly, we have $\mu=\mathbb{B}_t(\mu_t)$ if $t\geq0$. In addition, if
$0<t\leq\mathrm{Ind}(\mu)$ then the measure $\mu^{\uplus t}$ is
$\boxplus$-infinitely divisible and its Voiculescu transform can be
expressed as
\begin{equation}\label{N2}
\phi_{\mu^{\uplus t}}=E_{\mu_t^{\uplus t}}=
E_{(\mu^{\uplus(2t)})^{\boxplus1/2}}.
\end{equation}
\end{prop}

\begin{pf} It is clear that the measure $\mu_t$ in (\ref{N1}) is defined and (\ref{N}) holds if
$-1<t\leq0$. If $t=\mathrm{Ind}(\mu)<\infty$ then
$1-\mathrm{Ind}(\mu^{\uplus(t+1)})=1/(t+1)$ by Lemma \ref{3.4},
whence the measure $\mu_t$ in (\ref{N1}) is defined by Proposition
\ref{3.5}. The same result also holds for
$0<t<\mathrm{Ind}(\mu)\in(0,\infty]$, and hence the identity
$\mu=\mathbb{B}_t(\mu_t)$ holds. Next, observe that
$\mathrm{Ind}(\mu^{\uplus t})\geq1$ for $0<t\leq\mathrm{Ind}(\mu)$,
and therefore $\mu^{\uplus t}$ is $\boxplus$-infinitely divisible.
The first equality in (\ref{N2}) follows by replacing $\mu$ with
$\mu^{\uplus(t+1)}$ and letting $p=1/(t+1)$ in Proposition
\ref{3.5}, while the second equality follows from (\ref{bijection}).
\end{pf} \qed

Next, we relate the limit laws to $\boxplus$-divisibility
indicators.

\begin{prop} \label{limit} Let $q>0$ and $\{\mu_n\}$ be a sequence of measures in
$\mathcal{M}$ such that $\mu_n\to\mu$ weakly as $n\to\infty$ for
some $\mu\in\mathcal{M}$. Then the following statements hold.
\begin{enumerate} [$\qquad(1)$]
\item {The inequality
$\lim\sup_n\mathrm{Ind}(\mu_n)\leq\mathrm{Ind}(\mu)$ holds.}
\item {For any $p>0$ with
$1-\inf_n\mathrm{Ind}(\mu_n)\leq p$, $\mu_n^{\boxplus
p}\to\mu^{\boxplus p}$ weakly as $n\to\infty$.}
\item {The measures $\mu_n^{\uplus q}\to\mu^{\uplus q}$ weakly as
$n\to\infty$.}
\end{enumerate}
\end{prop}

\begin{pf} The measure $\mu_n^{\boxplus p}$ in (2) is defined for all
$n$ by Proposition \ref{3.5}, whence (2) holds by [\ref{HV2},
Proposition 5.7]. The assertion (3) holds by [\ref{BP}, Proposition
6.2]. To prove (1), first consider the case that
$0<t:=\lim\sup_n\mathrm{Ind}(\mu_n)<\infty$. Then for sufficiently
small $\epsilon>0$ we have
$1<\lim\sup_n\mathrm{Ind}(\mu_n^{\uplus(t-\epsilon)})$ by Lemma
\ref{3.4}, whence there exists a subsequence $\{\mu_{n_k}\}$ such
that $\mu_{n_k}^{\uplus(t-\epsilon)}$ is $\boxplus$-infinitely
divisible for all $k$. Since the set of $\boxplus$-infinitely
divisible measures is weakly closed, we see that
$\mathrm{Ind}(\mu^{\uplus(t-\epsilon)})\geq1$ by (3), which yields
$t-\epsilon\leq\mathrm{Ind}(\mu)$. Letting $\epsilon\to0$ shows
$t\leq\mathrm{Ind}(\mu)$. If $t=\infty$ then by similar arguments it
is easy to see that $m\leq\mathrm{Ind}(\mu)$ for any $m>0$, and
therefore $\mathrm{Ind}(\mu)=\infty$. The assertion (1) clearly
holds if $t=0$, and hence the proof is complete.
\end{pf} \qed

It was shown in Proposition \ref{3.1} that the subordination
function for the $\boxplus$-convolution power appearing in
(\ref{p-power}) is in fact the reciprocal Cauchy transform of some
$\boxplus$-infinitely divisible measure. The following theorem
states that the converse is also true. For other related results
about the $\boxplus$-infinite divisibility of the subordination
functions, we refer the reader to [\ref{G4}] and [\ref{Nica}].

\begin{thm} \label{sub} If $\mu\in\mathcal{M}$ then the following statements $(1)$ and $(2)$ are
equivalent.
\begin{enumerate} [$\qquad(1)$]
\item {The measure $\mu$ is $\boxplus$-infinitely divisible.}
\item {The function $F_\mu$ is the right inverse of some function
$H\in\mathcal{H}$.}
\end{enumerate}
If $(1)$ and $(2)$ hold then $F_\mu$ extends continuously to
$\mathbb{C}^+\cup\mathbb{R}$, $H$ can be written as
\begin{equation} \label{H}
H(z)=F_\mu^{-1}(z)=pz+(1-p)F_{\mu_p}(z),\;\;\;\;\;z\in\mathbb{C}^+,
\end{equation}
$\phi_\mu=E_{\mathbb{B}^{-1}(\mu)}$, and
\[F_{\mu^{\uplus
p^*}}(z)=F_{\mu_p}(F_\mu(z)),\;\;\;\;\;z\in\mathbb{C}^+,\] where
\begin{equation} \label{mup} \mu_p=\left(\mu^{\uplus
p^*}\right)^{\boxplus\frac{1}{p}},\;\;\;\;\;p>1.
\end{equation}
Moreover, for $r>0$ the measure $\mu^{\boxplus r}$ is
$\boxplus$-infinitely divisible and
\[\phi_{\mu^{\boxplus
r}}=E_{\left(\mu^{\uplus(1+r)}\right)^{\boxplus\frac{r}{1+r}}}.\]
\end{thm}

\begin{pf} First suppose that (2) holds, i.e., there exists some
$\nu\in\mathcal{M}$ such that $H(z)=2z-F_\nu(z)$ and
$H(F_\mu(z))=z$, $z\in\mathbb{C}^+$. By Proposition \ref{prop2.1}(5)
we see that $F_\mu(z)=\left[F_{\nu^{\boxplus2}}(z)+z\right]/2$ or,
equivalently, $\mu=\mathbb{B}_1(\nu)$, whence (1) holds by
Proposition \ref{3.1}. Conversely, if $\mu$ is $\boxplus$-infinitely
divisible then the measure $\mu_p$ in (\ref{mup}) is defined by
Lemma \ref{3.4} and Proposition \ref{3.5}. Moreover, by the fact
$\mathbb{B}_{p,1/p^*}(\mu_p)=\mu$ and Proposition \ref{3.1}, we
obtain $\phi_\mu=E_{\mu_p^{\uplus(p-1)}}$, which yields the
implications that (1) implies (2), and (\ref{H}). The last assertion
follows from (\ref{bijection}) and Proposition \ref{3.6}(1). Indeed,
we have
\[\phi_{\mu^{\boxplus r}}=rE_{\mathbb{B}^{-1}(\mu)}=E_{((\mu^{\uplus2})^{\boxplus1/2})^{\uplus r}}
=E_{(\mu^{\uplus(1+r)})^{\boxplus r/(1+r)}},\] as desired. This
finishes the proof.
\end{pf} \qed

Next, we analyze the supports and regularity of
$\boxplus$-infinitely divisible measures. Given such a measure
$\mu$, let $\Omega=F_\mu(\mathbb{C}^+)$. Then by Proposition
\ref{prop2.1}(1) and \ref{sub}, $\Omega=\{z\in\mathbb{C}^+:\Im
F_\mu^{-1}(z)>0\}$ is a simply connected domain and
$\partial\Omega=F_\mu(\mathbb{R})$ is the graph of the continuous
function
\begin{align*}f(x)&=\inf\left\{y>0:\Im F_\mu^{-1}(x+iy)>0\right\} \\
&=\inf\left\{y>0:\frac{\Im
E_{\mathbb{B}^{-1}(\mu)}(x+iy)}{y}>-1\right\},\;\;\;\;\;x\in\mathbb{R}.
\end{align*} Then Theorem \ref{Hthm}(1) shows that the function
$\psi(x)=F_\mu^{-1}(x+if(x))$, $x\in\mathbb{R}$, is homeomorphism on
$\mathbb{R}$. With the help of Proposition \ref{prop2.1}, we have
the following conclusions.

\begin{thm} \label{3.10} Suppose that $\mu\in\mathcal{M}$ is $\boxplus$-infinitely divisible.
\begin{enumerate} [$\qquad(1)$]
\item {The function $\phi_\mu$ has a continuous extension to
$\partial\Omega$ and for any $z_1,z_2\in\overline{\Omega}$,
\[|\phi_\mu(z_1)-\phi_\mu(z_2)|\leq|z_1-z_2|.\]}
\item {For any $z_1,z_2\in\mathbb{C}^+\cup\mathbb{R}$,
\[\frac{|z_1-z_2|}{2}\leq|F_\mu(z_1)-F_\mu(z_2)|.\]
Consequently, $G_\mu$ has a continuous extension to $\mathbb{R}$
except one point and the measure $\mu$ has at most one atom.}
\item {The measure $\mu$ has an atom if and only if
$0\in\partial\Omega$ and
\[\lim_{\epsilon\downarrow0}\frac{F_\mu^{-1}(i\epsilon)-F_\mu^{-1}(0)}{i\epsilon}=m>0,\] in
which case the point $F_\mu^{-1}(0)$ is an atom of $\mu$ with mass
$m$.}
\item {The nonatomic part of $\mu$ is absolutely continuous $($with respect to Lebesgue
measure$)$.}
\item {The measure $\mu^{\mathrm{ac}}$ is concentrated on the
set $\psi(\overline{V^+})$, where $V^+=\{x:f(x)>0\}$.}
\item {At the point $\psi(x)$, $x\in V^+$, the density of
$\mu^{\mathrm{ac}}$ is analytic and given by
\[\frac{d\mu^{\mathrm{ac}}}{dx}(\psi(x))=\frac{f(x)}{\pi(x^2+f^2(x))}.\]}
\item {The measure $\mu$ is compactly supported if and only if so is
$f$.}
\end{enumerate}
\end{thm}

\begin{pf} By letting $p=2$ in Proposition \ref{prop2.1} and Theorem \ref{sub},
we have $H(z)=F_\mu^{-1}(z)=2z-F_{\mathbb{B}^{-1}(\mu)}(z)$,
$z\in\mathbb{C}^+$. Then it follows from \ref{prop2.1}(2) that
\begin{equation} \label{in}
\int_\mathbb{R}\frac{s^2+1}{|s-z|^2}\;d\sigma(s)\leq1,\;\;\;\;\;z\in\overline{\Omega}.
\end{equation} where $\sigma$ is the measure in the
Nevanlinna representation of $F_{\mathbb{B}^{-1}(\mu)}$. Since
$\phi_\mu=E_{\mathbb{B}^{-1}(\mu)}$, the inequality in (1) holds for
$z\in\Omega$ by H\"{o}lder inequality and (\ref{in}), whence (1)
holds by continuous extension. The assertion (2) follows from
Proposition \ref{prop2.1} (3). Observer that $\mu$ has an atom at
$\alpha$ if and only if $F_\mu(\alpha)=0$ and the
Julia-Carath\'{e}odory derivative $F_\mu'(\alpha)<\infty$, which
happens if and only if $0\in\partial\Omega$ and
\[0<\frac{1}{F_\mu'(\alpha)}=(F_\mu^{-1})'(0),\] where $(F_\mu^{-1})'(0)$
is the Julia-Carath\'{e}odory derivative of $F_\mu^{-1}$ at $0$.
Hence $\mu(\{\alpha\})=(F_\mu^{-1})'(0)$ and (3) holds. Next, note
that for any $x\in\mathbb{R}$ we have $F_\mu(\psi(x))=x+if(x)$.
Since $F_\mu$ extends continuously to $\mathbb{C}^+\cup\mathbb{R}$,
applying the inversion formula (\ref{inversion}) gives
\[\frac{d\mu^{\mathrm{ac}}}{dx}(\psi(x))=\frac{-1}{\pi}\Im G_\mu(\psi(x))=
\frac{f(x)}{\pi(x^2+f^2(x))},\;\;\;\;\;x\in V^+,\] which, along with
\ref{prop2.1}(4) gives (5) and (6). As noted above, $F_\mu(x)=0$
a.e. relative to the singular part of $\mu$, from which we deduce
that the singular part of $\mu$ is atomic, which gives (4). That (7)
follows from (5) and the fact that $\mu$ has at most one atom.
\end{pf} \qed

The constants appearing in \ref{3.10}(1) and (2) are sharp. Indeed,
by considering the standard semicircular distribution we have
$\Omega=\{z\in\mathbb{C}^+:|z|>1\}$, and then taking $z_1=1$ and
$z_2=-1$ shows that $1$ is the best constant in (1), whence the same
conclusion for (2) follows immediately.

Recall that the compound free Poisson distribution $p(\lambda,\nu)$
with the rate $\lambda>0$ and jump distribution $\nu$is defined as
the weak limit as $n\to\infty$ of $\mu_n^{\boxplus n}$, where
\[\mu_n=\left(1-\frac{\lambda}{n}\right)\delta_0+\frac{\lambda}{n}\nu\] and $\nu$ is compactly
supported. The next proposition generalizes the jump distribution
with compact support to any measure in $\mathcal{M}$.

\begin{prop} \label{Poisson} Given $\nu\in\mathcal{M}$, define
\[d\rho(s)=\frac{s^2}{s^2+1}\;d\nu(s).\] Then $p(\lambda,\nu)=\mathbb{B}_{1+\lambda,1/(1+\lambda)^*}(\mu_0)$, where
$\mu_0$ is a measure in $\mathcal{M}$ whose reciprocal Cauchy
transform satisfies
\[F_{\mu_0}(z)=-\int_\mathbb{R}\frac{s}{s^2+1}\;d\nu(s)+z+\int_\mathbb{R}\frac{1+sz}{s-z}\;d\rho(s).\]
Consequently, $p(\lambda,\nu)$ is a $\boxplus$-infinitely divisible
measure with an atom at $0$ of mass $1-\lambda$ for $\lambda<1$ and
no atom for $\lambda\geq1$,
\[\phi_{p(\lambda,\nu)}(z)=\lambda E_{\mu_0}(z)=\lambda z\int_\mathbb{R}\frac{s}{z-s}\;d\nu(s),\] and
\[\mathrm{Ind}(p(\lambda,\nu))=\frac{\mathrm{Ind}(\mu_0)+\lambda}{\lambda}.\]
\end{prop}

\begin{pf} Since $s^2/(s^2+1)\in L^1(d\nu)$, the measure $\rho$ is finite and positive, and
the limit
\[\lim_{n\to\infty}\int_\mathbb{R}\frac{ns}{s^2+1}\;d\mu_n(s)=\int_\mathbb{R}\frac{\lambda s}{s^2+1}\;d\nu(s)\]
exists. Moreover, it is easy to see that
\[\frac{ns^2}{s^2+1}\;d\mu_n(s)\to \lambda d\rho(s)\] weakly. By Theorem
\ref{thm2.2}, the measure $\mu_n^{\boxplus n}$ converges weakly to
$p(\lambda,\nu)$, which satisfies
\[\phi_{p(\lambda,\nu)}(z)=\int_\mathbb{R}\frac{\lambda s}{s^2+1}\;d\nu(s)+
\int_\mathbb{R}\frac{\lambda(1+sz)}{z-s}\;d\rho(s)=\lambda
z\int_\mathbb{R}\frac{s}{z-s}\;d\nu(s).\] On the other hand, the
definition of $\mu_0$ and Proposition 3.1 show that
\[\phi_{\mathbb{B}_{1+\lambda,1/(1+\lambda)^*}
(\mu_0)}=\lambda E_{\mu_0}=\phi_{p(\lambda,\nu)}.\] Then by Lemma
\ref{3.4} and (\ref{Nindicator}) we have
\[\mathrm{Ind}(p(\lambda,\nu))=\mathrm{Ind}(\mathbb{B}_\lambda(\mu_0)^{\uplus \lambda})
=\frac{\mathrm{Ind}(\mu_0)+\lambda}{\lambda}.\]
Next, we apply Theorem \ref{3.10}(3) to locate the atom of
$p(\lambda,\nu)$. Since $\phi_{p(\lambda,\nu)}(0)=0$, $0\in\partial
F_{p(\lambda,\nu)}(\mathbb{C}^+)$. Moreover, by the dominated
convergence theorem we obtain
\[\lim_{\epsilon\downarrow0}\frac{\phi_{p(\lambda,\nu)}(i\epsilon)-\phi_{p(\lambda,\nu)}(0)}{i\epsilon}
=\lambda\lim_{\epsilon\downarrow0}\int_\mathbb{R}\frac{s}{i\epsilon-s}\;d\nu(s)=-\lambda,\]
which gives the desired result. This completes the proof.
\end{pf} \qed

Since the $\boxplus$-divisibility indicator is zero for any measure
with finite support, we have the following result.

\begin{cor} We have
$p(\lambda,\delta_a)=\mathbb{B}_{1+\lambda,1/(1+\lambda)^*}(\mu_0)$,
where $\mu_0=(\delta_0+\delta_{2a})/2$,
$\mathrm{Ind}(p(\lambda,\delta_a))=1$, and
$\phi_{p(\lambda,\delta_a)}(z)=a\lambda z/(z-a)$.
\end{cor}

We finish this section with an interesting observation. If
$\mathrm{Ind}(\mu)>1$ then
$\mathrm{Ind}(\mathbb{B}^{-1}(\mu))=\mathrm{Ind}(\mu)-1>0$ by
(\ref{Nindicator}) and (\ref{bijection}). This implies that
$\phi_\mu=E_{\mathbb{B}^{-1}(\mu)}$ has a continuous extension to
$\mathbb{C}^+\cup\mathbb{R}$ by Proposition \ref{3.5}, whence we
have the following proposition.

\begin{prop} If $\mu\in\mathcal{M}$ with $\mathrm{Ind}(\mu)>1$ then
$\phi_\mu$ has a continuous extension to
$\mathbb{C}^+\cup\mathbb{R}$.
\end{prop}

\setcounter{equation}{0}
\section{Measures with mean zero and finite variance}

Recall that for $\mu\in\mathcal{M}$ with mean zero and unit
variance, $\Phi(\mu)$ is the unique measure in $\mathcal{M}$
satisfying the Eq. (\ref{EMaa}) with $\sigma^2=1$, i.e.,$
E_\mu=G_{\Phi(\mu)}$. In general, a measure $\mu$ has mean $m$ and
finite variance $\sigma^2$ if and only if $\mu\boxplus\delta_{-m}$
has mean zero and variance $\sigma^2$ because
$d(\mu\boxplus\delta_{-m})(s)=d\mu(s+m)$, and hence
$E_{\mu\boxplus\delta_{-m}}=\sigma^2G_{\Phi((\mu\boxplus\delta_{-m})^{\uplus1/\sigma^2})}$.
Since $\mathrm{Ind}(\mu)=\mathrm{Ind}(\mu\boxplus\delta_a)$ for any
$a\in\mathbb{R}$ by [\ref{Japan}, Proposition 3.7], in what follows
we only consider measures with mean zero and finite variance.

Recall that the free Brownian motion started at $\nu\in\mathcal{M}$
is the process $\{\nu\boxplus\gamma_t:t\geq0\}$. The connection
among this process, the map $\mathbb{B}_t$, and the subordination
function of the $\boxplus$-convolution powers is described in the
following theorem, which was proved in [\ref{BN1}] and
[\ref{Biane1}]. For the completeness, we provide its statement and
proof.

\begin{thm} \label{Brownian} If $\mu\in\mathcal{M}$ has mean zero
and variance $\sigma^2$, and $\nu=\Phi(\mu^{\uplus1/\sigma^2})$ then
\begin{equation} \label{Gmotion}
G_{\nu\boxplus\gamma_{t\sigma^2}}(z)=
G_\nu(F_{\mathbb{B}_{t+1,t/(t+1)}(\mu)}(z))=\frac{E_{\mathbb{B}_{t+1,t/(t+1)}(\mu)}(z)}{t\sigma^2},\;\;\;\;\;
z\in\mathbb{C}^+\cup\mathbb{R},
\end{equation} where $t>0$. Consequently, we have
$\phi_{\mathbb{B}_{t+1,t/(t+1)}(\mu)}=t\sigma^2G_\nu$ and
\begin{equation} \label{motion}
E_{\mathbb{B}_t(\mu)}=\sigma^2G_{\nu\boxplus\gamma_{t\sigma^2}}.
\end{equation}
\end{thm}

\begin{pf} Let $p=t+1>1$. Since $E_\mu=\sigma^2G_\nu$, it follows that
\[H_p(z):=pz+(1-p)F_\mu(z)=z+(p-1)\sigma^2G_\nu(z),\;\;\;\;\;z\in\mathbb{C}^+.\]
If $\omega_p$ is the right inverse of $H_p$ then [\ref{Biane1},
Proposition 2] shows that
\[G_{\nu\boxplus\gamma_{(p-1)\sigma^2}}(z)=G_\nu(\omega_p(z))=\frac{z-\omega_p(z)}{(p-1)\sigma^2},
\;\;\;\;\;z\in\mathbb{C}^+\cup\mathbb{R}.\] Since
$\omega_p=F_{\mathbb{B}_{p,1/p^*}(\mu)}$ by Proposition \ref{3.1},
the above identity yields (\ref{Gmotion}). Finally, the rest
assertions follow from
$\phi_{\mathbb{B}_{p,1/p^*}}=E_{\mu^{\uplus(p-1)}}=(p-1)\sigma^2G_\nu$
and (\ref{Gmotion}).
\end{pf} \qed

The identity (\ref{motion}) indicates that $\mu^{\boxplus p}$ has
mean zero and finite variance $p\sigma^2$ if $p\geq1$. The next
result shows that this is also true for the measure $\mu^{\boxplus
p}$ whenever it is defined.

\begin{lem} \label{4.1} Suppose that $\mu\in\mathcal{M}$ has mean zero and variance $\sigma^2$.
If the measure $\mu^{\boxplus p}$ is defined for some $p>0$ then it
has mean zero and variance $p$, in which case
\begin{align*}
&\Phi\left(\left(\mu^{\boxplus
p}\right)^{\uplus1/(p\sigma^2)}\right)=\Phi\left(\mu^{\uplus1/\sigma^2}\right)
\boxplus\gamma_{(p-1)\sigma^2},\;\;\;\;\;p\geq1, \\
&\Phi\left(\left(\mu^{\boxplus
p}\right)^{\uplus1/(p\sigma^2)}\right)\boxplus\gamma_{(1-p)\sigma^2}=\Phi
\left(\mu^{\uplus1/\sigma^2}\right),\;\;\;\;\;p<1.
\end{align*}
\end{lem}

\begin{pf} By (\ref{motion}), it suffices to show the lemma for the case
$1-\mathrm{Ind}(\mu)\leq p<1$. Let
$\nu=\Phi(\mu^{\uplus1/\sigma^2})$, $\mu_p=\mu^{\boxplus p}$, and
$H(z)=z/p+(1-1/p)F_{\mu_p}(z)$. Then it follows from (\ref{p<1})
that $F_{\mu_p}(iy)=F_\mu(H(iy))$ or, equivalently,
$E_{\mu_p}(iy)=iy-H(iy)+\sigma^2G_\nu(H(iy))$ for sufficiently large
$y>0$. Since $z-H(z)=(1-1/p)E_\mu(z)$, we see that
$E_{\mu_p}(iy)=p\sigma^2G_\nu(H(iy))$ for sufficiently large $y>0$.
Next, we claim that $E_{\mu_p}/(p\sigma^2)\in\mathcal{G}$. Indeed,
since $\lim_{y\to\infty}H(iy)/(iy)=1$, for any $\alpha>0$ there
exists a number $\beta>0$ such that
\[\left|\frac{H(iy)}{iy}-1\right|<c_\alpha:=\frac{\alpha}{\sqrt{1+\alpha^2}},\;\;\;\;\;y>\beta,\]
from which we deduce that $H(iy)\in\Gamma_{\alpha,\beta'}$ for
$y>\beta$, where $\beta'=(1-c_\alpha)\beta$. By [\ref{HV2},
Proposition 5.1], we obtain
\[\lim_{y\to\infty}iyE_{\mu_p}(iy)=p\sigma^2\left(\lim_{y\to\infty}\frac{iy}{H(iy)}\right)
\left(\lim_{y\to\infty}H(iy)G_\nu(H(iy))\right)=p\sigma^2,\] which
yields that $E_{\mu_p}/(p\sigma^2)\in\mathcal{G}$, as desired.
Finally, let $\nu_p=\Phi((\mu^{\boxplus p})^{\uplus1/(p\sigma^2)})$.
Then $E_{\mu^{\boxplus p}}=p\sigma^2G_{\nu_p}$ and (\ref{motion})
show that
\[p\sigma^2G_{\Phi(\mu^{\uplus1/\sigma^2})}=E_{\mu^{\uplus
p}}=E_{\mathbb{B}_{1/p-1}(\mu^{\boxplus
p})}=p\sigma^2G_{\nu_p\boxplus\gamma_{(1-p)\sigma^2}},\] which gives
the last assertion. If $p>1$ then
$E_\mu=\sigma^2G_{\Phi(\mu^{\uplus1/\sigma^2})}$ and (\ref{motion})
yield
\[E_{\mathbb{B}_{p-1}(\mu)}=\sigma^2G_{\Phi(\mu^{\uplus1/\sigma^2})\boxplus\gamma_{(p-1)\sigma^2}},\]
as desired. This completes the proof.
\end{pf} \qed

\begin{prop} \label{meanbasic} Suppose that $\mu\in\mathcal{M}$ has mean zero and
finite variance $\sigma^2$. If $t$ is a finite number with $0\leq
t\leq \varphi(\mu)$ and $\mu_t$ is the measure defined in
$(\ref{N1})$ then
\[E_{\mu_t}=\sigma^2G_{\nu_t}\;\;\;\;\;\mathrm{and}\;\;\;\;\;
E_\mu=\sigma^2G_{\nu_t\boxplus\gamma_{t\sigma^2}},\] where
$\nu_t=\Phi\left(\mu_t^{\uplus1/\sigma^2}\right)$.
\end{prop}

\begin{pf} By Lemma \ref{4.1}, it is clear that $\mu_t$ has mean
zero and variance $\sigma^2$, whence the conclusions follows from
(\ref{motion}).
\end{pf} \qed

The preceding proposition gives a reformulation for the
$\boxplus$-divisibility indicator of measures with mean zero and
finite variance.

\begin{cor} \label{Ind} If $\mu\in\mathcal{M}$ has mean zero and finite variance
$\sigma^2$ then
\[\mathrm{Ind}(\mu)=\sup\left\{t\geq0:E_\mu=\sigma^2G_{\nu_t\boxplus\gamma_{t\sigma^2}}\;\;\mathrm{for}\;\;
\mathrm{some}\;\;\nu_t\in\mathcal{M}\right\}.\]
\end{cor}

The preceding corollary enables us to associate to each measure
$\nu\in\mathcal{M}$ a nonnegative number:
\[C(\nu)=\sup\{t\geq0:\nu=\nu_t\boxplus\gamma_t\;\;\mathrm{for\;\;some\;\;}\nu_t\in\mathcal{M}\}.\]
We will call $C(\nu)$ the semicircular decomposition indicator of
$\nu$. The connection between $\boxplus$-divisibility indicator and
semicircular decomposition indicator is described in the next
result.

\begin{thm} \label{BInd} For any $\nu\in\mathcal{M}$ we have
$B(\nu)=\mathrm{Ind}(\Phi^{-1}(\nu))$.
\end{thm}

We now characterize $\boxplus$-infinitely divisible measures with
mean zero and finite variance.

\begin{thm} \label{divisible} If $\mu\in\mathcal{M}$ and $\sigma\in(0,\infty)$ then the
following statements are equivalent:
\begin{enumerate} [$\qquad(1)$]
\item {$\mu$ is a $\boxplus$-infinitely divisible measure with mean
zero and variance $\sigma^2$;}
\item {there exists a measure $\nu\in\mathcal{M}$ such that
$\phi_\mu=\sigma^2G_\nu$;}
\item {$F_\mu$ is the right inverse of
some $H\in\mathcal{H}$ satisfying
$\lim_{y\to\infty}iy(H(iy)-iy)=\sigma^2$;}
\item {there exists a measure $\nu\in\mathcal{M}$ such that
$E_\mu=\sigma^2G_{\nu\boxplus\gamma_{\sigma^2}}$.}
\end{enumerate} If $(1)$-$(4)$ hold and $p=1+\sigma^2$ then the
measure $\nu$ in $(2)$ and $(4)$ can be expressed as
\[\nu=\Phi\left(\left(\mu^{\uplus p^*}\right)^{\boxplus1/p}\right).\] The function $H$ in $(3)$ can be expressed
as \begin{equation} \label{HG} H(z)=z+\sigma^2G_\nu(z),
\end{equation} and
\begin{equation} \label{HG2} G_{\nu\boxplus\gamma_{\sigma^2}}(z)
=G_\nu(F_\mu(z)),\;\;\;\;\;z\in\mathbb{C}^+\cup\mathbb{R}.
\end{equation} Moreover, for any $r>0$ we have
\[E_{\mu^{\boxplus r}}=r\sigma^2G_{\nu\boxplus\gamma_{r\sigma^2}}.\]
\end{thm}

\begin{pf} First suppose that (1) holds. Then the measure $(\mu^{\uplus2})^{\boxplus1/2}=\mathbb{B}^{-1}(\mu)$
has mean zero and variance $\sigma^2$ by Lemma \ref{4.1}, whence
$\phi_\mu=E_{\mathbb{B}^{-1}(\mu)}=\sigma^2G_\nu$ for some
$\nu\in\mathcal{M}$ and (2) follows. The definition of $\Phi$ shows
that $\nu$ can be expressed as
\[\nu=\Phi\left((\mathbb{B}^{-1}(\mu))^{\uplus1/\sigma^2}\right)=\Phi\left(\left(\mu^{\uplus q^*}\right)^
{\boxplus1/q}\right),\] where the Eq. (\ref{formula2}) is used in
the second equality above. If (2) holds then
$H(z)=F_\mu^{-1}(z)=\phi_\mu(z)+z$, which implies (3). If the
statement (3) holds then $H(z)=z+\sigma^2G_{\nu_1}$ for some
$\nu_1\in\mathcal{M}$. Then [\ref{Biane1}, Proposition 2] shows that
$F_\mu$ is the subordination function of
$\nu\boxplus\gamma_{\sigma^2}$ with respect to $\nu$, whence we have
\[G_{\nu_1\boxplus\gamma_{\sigma^2}}(z)=G_{\nu_1}(F_\mu(z))=\frac{z-F_\mu(z)}{\sigma^2},\;\;\;\;\;z\in\mathbb{C}^+,\]
and the assertion (4) holds. The implication that (4) implies (1)
follows from Corollary \ref{Ind}. Moreover, the identity
(\ref{motion}) shows that
$E_\mu=E_{\mathbb{B}(\mathbb{B}^{-1}(\mu))}=\sigma^2G_{\nu_2\boxplus\gamma_{\sigma^2}}$,
whence the assertions (\ref{HG}) and (\ref{HG2}) hold by the
preceding discussions. For the last assertion it suffices to show
that $\nu_r:=\Phi((\mu^{\boxplus r})^{\uplus1/(r\sigma^2)})=
\nu\boxplus\gamma_{r\sigma^2}$. If $r<1$ then
$\nu_r\boxplus\gamma_{(1-r)\sigma^2}=\Phi(\mu^{\uplus1/\sigma^2})$
by Lemma \ref{4.1}. Since
$\Phi(\mu^{\uplus1/\sigma^2})=\nu\boxplus\gamma_{\sigma^2}$, the
desired equality follows. Similarly, if $r>1$ then
$\nu_r=\Phi(\mu^{\uplus1/\sigma^2})\boxplus\gamma_{(r-1)\sigma^2}=
\nu\boxplus\gamma_{r\sigma^2}$, as desired.
\end{pf} \qed

Let $H$ be the function defined as in (\ref{HG}) and
\[\Omega=\{z\in\mathbb{C}^+:\Im H(z)>0\}.\] It was shown in [\ref{Biane1}] that the function $G_\nu$ extends
continuously to $\overline{\Omega}$ and this extension is Lipschitz
continuous on $\overline{\Omega}$ with the Lipschitz constant
$1/\sigma^2$. Moreover,
\[|G_\nu(z)|\leq\frac{1}{\sigma},\;\;\;\;z\in\overline{\Omega}.\]
Combining these facts and Theorem \ref{divisible} gives the
following result.

\begin{cor} If $\mu\in\mathcal{M}$ is a $\boxplus$-infinitely
divisible measure with mean zero and finite variance $\sigma^2$ then
\[|G_{\nu\boxplus\gamma_{\sigma^2}}(z_1)-G_{\nu\boxplus\gamma_{\sigma^2}}(z_2)|\leq
\frac{1}{\sigma^2}|F_\mu(z_1)-F_\mu(z_2)|,\;\;\;\;\;z_1,z_2\in\mathbb{C}^+\cup\mathbb{R},\]
and
\[|\phi_\mu(z)|\leq\sigma,\;\;\;\;\;
z\in\overline{\Omega},\] where $\nu=\Phi(\mu^{\uplus1/\sigma^2})$
and $\Omega=F_\mu(\mathbb{C}^+)$.
\end{cor}

It was shown before that $E_\mu$ has a continuous extension to
$\mathbb{C}^+\cup\mathbb{R}$ if $\mathrm{Ind}(\mu)>0$. In general,
the converse is not true. Indeed, let $\mu\in\mathcal{M}$ be so that
$E_\mu=G_N=1/(z+i)$, where $N$ is the Cauchy distribution. Since
$\phi_N=-i$, it is easy to see that $N$ cannot be written as a free
Brownian motion stated at some measure, whence $B(N)=0$, which
yields $\mathrm{Ind}(\mu)=0$ by Theorem \ref{BInd}. In the following
theorem, we improve this result for measures with mean zero and
finite variance.

\begin{thm} \label{NandS} If $\mu\in\mathcal{M}$ has mean zero and
finite variance $\sigma^2$ then
\begin{enumerate} [$\qquad(1)$]
\item {$\mathrm{Ind}(\mu)>0$ if and only if $E_\mu=\sigma^2G_{\nu\boxplus\gamma_t}$ for some $\nu\in\mathcal{M}$
and $t>0$;}
\item {$\mathrm{Ind}(\mu)>1$ if and only if
$\phi_\mu=\sigma^2G_{\nu\boxplus\gamma_t}$ for some
$\nu\in\mathcal{M}$ and $t>0$}.
\end{enumerate}
\end{thm}

\begin{pf} The assertion (1) was proved in Proposition
\ref{meanbasic}. Since $\phi_\mu=E_{\mathbb{B}^{-1}(\mu)}$ and
$\mathrm{Ind}(\mu)=1+\mathrm{Ind}(\mathbb{B}^{-1}(\mu))$ if $\mu$ is
$\boxplus$-infinitely divisible, the assertion (2) follows (1).
\end{pf} \qed

For $a\in\mathbb{R}$, by the fact $(\mu\boxplus\delta_a)^{\boxplus
p}=\mu^{\boxplus p}\boxplus\delta_{pa}$ and the identity
$(\mu\boxplus\delta_a)^{\uplus q}=(\mu^{\uplus
q}\boxplus\delta_a)\uplus\delta_{(q-1)a}$ shown in [\ref{Japan},
Proposition 3.7] we have
\begin{equation} \label{Bpq} \mathbb{B}_{p,q}(\mu\boxplus\delta_a)=
(\mathbb{B}_{p,q}(\mu)\boxplus\delta_{pa})\uplus\delta_{p(q-1)a}.
\end{equation} Next, we use (\ref{Bpq}) to investigate the free compound
Poisson distribution $p(\lambda,\nu)$, where $\nu$ has finite
variance.

\begin{prop} Suppose that $\nu\in\mathcal{M}$ has mean $m$ and finite variance
$\sigma^2$. Then
\[E_{p(\lambda,\nu)\boxplus\delta_{-\lambda m}}=
\lambda m_2G_{\nu_0\boxplus\gamma_{\lambda m_2}}\] and
\[\phi_{p(\lambda,\nu)\boxplus\delta_{-\lambda
m}}=\lambda m_2G_{\nu_0},\] where $m_2=m^2+\sigma^2$ is the second
moment of $\nu$ and $d\nu_0(s)=s^2/m_2d\nu(s)$. Consequently,
$p(\lambda,\nu)$ has mean $\lambda m$ and variance $\lambda m_2$,
and $\mathrm{Ind}(p(\lambda,\nu))>1$ if $B(\nu_0)>0$.
\end{prop}

\begin{pf} If $\mu_0$ is the measure defined in Proposition \ref{Poisson} then
\begin{equation} \label{4.2} E_{\mu_0}(z)=\int_\mathbb{R}s\;d\nu(s)+
\int_\mathbb{R}\frac{s^2}{z-s}\;d\nu(s)=m+m_2G_{\nu_0}(z),
\end{equation} from which we obtain
\begin{align*}
E_{\mu_0\boxplus\delta_{-m}}(z)&=E_{\mu_0}(z+m)-m \\
&=m_2G_{\nu_0}(z+m)=m_2G_{\nu_0\boxplus\delta_{-m}}(z).
\end{align*}
Then Theorem \ref{Brownian} shows that for any $p>1$ we have
\begin{equation} \label{4.3} E_{\mathbb{B}_{p,1/p^*(\mu_0\boxplus\delta_{-m})}}
=(p-1)E_{\mathbb{B}_{p-1}(\mu_0\boxplus\delta_{-m})}
=(p-1)m_2G_{\nu_0\boxplus\delta_{-m}\boxplus\gamma_{(p-1)m_2}}.
\end{equation} On the other hand, by (\ref{Bpq}) we have
\[E_{\mathbb{B}_{p,1/p^*(\mu_0\boxplus\delta_{-m})}}(z)
=E_{\mathbb{B}_{p,1/p^*}(\mu_0)}(z-pm)+(1-p)m,\] from which, along
with (\ref{4.3}), we deduce that
\[E_{\mathbb{B}_{p,1/p^*}(\mu_0)}(z)+(1-p)m=(p-1)m_2G_{\nu_0\boxplus\delta_{-m}\boxplus\gamma_{(p-1)m_2}}(z-pm)\]
or, equivalently,
\[E_{\mathbb{B}_{p,1/p^*}(\mu_0)\boxplus\delta_{(1-p)m}}=
(p-1)m_2G_{\nu_0\boxplus\delta_{(p-1)m}\boxplus\gamma_{(p-1)m_2}}.\]
Letting $p=\lambda+1$ in the above identity gives that
$p(\lambda,\nu)$ have mean $\lambda m$ and variance $\lambda m_2$.
Since $\phi_{p(\lambda,\nu)}=\lambda E_{\mu_0}$, it follows from
(\ref{4.2}) that
\[\phi_{p(\lambda,\nu)\boxplus\delta_{-\lambda
m}}=\lambda m_2G_{\nu_0}.\] The last assertion follows from
[\ref{Japan}, Proposition 3.7] and Corollary \ref{NandS}.
\end{pf} \qed

From the preceding proposition, it is easy to see that
$p(\lambda,\delta_a)$ has mean $\lambda a$ and variance $\lambda
a^2$, and $\mathrm{Ind}(p(\lambda,\delta_a))=1$ since
$\nu_0=\delta_a$.

\setcounter{equation}{0}
\section{Support and regularity for measures in $\mathbb{B}_{p,q}(\mathcal{M})$}

If $p,q>0$ then the measure $\mathbb{B}_{p,q}(\mu)$ (if
$\mu^{\boxplus p}$ is defined) is a Dirac measure $\delta_a$ if and
only if $\mu=\delta_{a/(pq)}$. For the rest of the paper we confine
our attention to the case of $\mu\in\mathcal{M}$ which is not a
point mass and follow the notations used in Proposition
\ref{prop2.1} and \ref{Hthm}. We denote by $\rho$ the unique nonzero
(because $\mu\neq\delta_a$) measure in the Nevanlinna representation
(\ref{NeF}) of $F_\mu$. Therefore, the Nevanlinna representation of
$F_{\mu^{\uplus q}}$ is
\begin{equation} \label{Neq}
F_{\mu^{\uplus q}}(z)=q\Re
F_\mu(1)+z+q\int_\mathbb{R}\frac{1+sz}{s-z}\;d\rho(s),\;\;\;\;\;z\in\mathbb{C}^+.
\end{equation}

In certain situation, $F_\mu$ is defined and takes a real value at
some $x\in\mathbb{R}$ (for instance, $x$ is an atom of $\mu$), in
which case we write $F_\mu(x)\in\mathbb{R}$. The following result
shows that for $p>1,q>0$, $F_{\mu^{\uplus q}}$ is Lipschitz
continuous on $\overline{\Omega_p}$ and takes real values on
$\overline{\Omega_p}\cap\mathbb{R}$.

\begin{prop} \label{Lip} For $p>1,q>0$, $F_{\mu^{\uplus
q}}$ extends continuously to $\overline{\Omega_p}$ and satisfies
\[\left|\frac{F_{\mu^{\uplus q}}(z_1)-F_{\mu^{\uplus q}}(z_2)}{z_1-z_2}\right|
\leq1+\frac{q}{p-1},\;\;\;\;z_1,z_2\in\overline{\Omega_p}.\]
Moreover, $(\ref{Neq})$ holds for $z\in\overline{\Omega_p}$ and the
Julia-Carath\'{e}odory derivative $F_{\mu^{\uplus q}}'$ is
\begin{equation} \label{JCq}
F_{\mu^{\uplus
q}}'(z)=1+q\int_\mathbb{R}\frac{s^2+1}{(s-z)^2}\;d\rho(s),\;\;\;\;\;z\in\overline{\Omega_p}.
\end{equation}
\end{prop}

\begin{pf} First, applying Proposition \ref{prop2.1}(2) and the H\"{o}lder inequality
to $E_\mu$ gives
\[\left|\frac{E_\mu(z_1)-E_\mu(z_2)}{z_1-z_2}\right|\leq\int_\mathbb{R}\frac{(s^2+1)d\rho(s)}{|s-z_1||s-z_2|}
\leq\frac{1}{p-1},\;\;\;\;\;z_1,z_2\in\Omega_p.\] Then by the
continuous extension, the above inequality holds for
$z_1,z_2\in\overline{\Omega_p}$, and therefore the Nevanlinna
representation (\ref{NeE}) of $E_\mu$ holds for
$z\in\overline{\Omega_p}$. Using the dominated convergence theorem,
the Julia-Carath\'{e}odory $E_\mu'$ is then given by
\[E_\mu'(z)=\lim_{\epsilon\downarrow0}
\frac{E_\mu(z+i\epsilon)-E_\mu(z)}{i\epsilon}
=-\int_\mathbb{R}\frac{s^2+1}{(s-z)^2}\;d\rho(s),\;\;\;\;\;z\in\overline{\Omega_p},\]
whence the desired results follow from the identities
$E_{\mu^{\uplus q}}=qE_\mu$ and $F_{\mu^{\uplus
q}}'=1-E_{\mu^{\uplus q}}'$.
\end{pf} \qed

The following lemma plays an important role in the investigation of
atoms of the measure $\mathbb{B}_{p,q}(\mu)$.

\begin{lem} \label{basic} Let $x\in\mathbb{R}$ and $f_\mu$ be the function defined as in $(\ref{fmu})$. Then
\\ $(1)$ $F_\mu(x)\in\mathbb{R}$ and the Julia-Carath\'{e}odory derivative
$F_\mu'(x)\in(1,\infty)$
\\ if and only if
\\ $(2)$ $F_\mu(x)\in\mathbb{R}$ and $f_\mu(x)\in(0,\infty)$,
\\ in which case $(\ref{NeF})$ holds for $z=x$ and $F_\mu'(x)=1+f_\mu(x)$.
\end{lem}

\begin{pf} First, suppose that (2) holds. Then
$x\in\overline{\Omega_p}$ for some $p>1$ and (\ref{NeF}) holds for
$z=x$ by Proposition \ref{prop2.1}(2). As shown in Proposition
\ref{Lip}, we have the Julia-Catath\'{e}odory $E_\mu'(x)=f_\mu(x)$,
whence Julia-Carath\'{e}odory derivative
$F_\mu'(x)=1+f_\mu(x)\in(1,\infty)$ and (1) follows. On the other
hand, if both $F_\mu(x)$ and the Julia-Carath\'{e}odory derivative
$F_\mu'(x)$ are real numbers then
\begin{align*}
F_\mu'(x)&=\lim_{\epsilon\downarrow0}\frac{\Re[F_\mu(x+i\epsilon)-F_\mu(x)]}{i\epsilon}
+\lim_{\epsilon\downarrow0}\frac{i\Im[F_\mu(x+i\epsilon)-F_\mu(x)]}{i\epsilon} \\
&=\lim_{\epsilon\downarrow0}\frac{\Im[F_\mu(x+i\epsilon)-F_\mu(x)]}{\epsilon}
=\lim_{\epsilon\downarrow0}\frac{\Im F_\mu(x+i\epsilon)}{\epsilon} \\
&=\lim_{\epsilon\downarrow0}\left(1+\int_{\mathbb{R}}\frac{s^2+1}{(s-x)^2+\epsilon^2}\;d\rho(s)\right) \\
&=1+\int_{\mathbb{R}}\frac{s^2+1}{(s-x)^2}\;d\rho(s),
\end{align*} where the monotone convergence theorem is used in the last equality.
This yields the implication that (1) implies (2) and the proof is
complete.
\end{pf} \qed

Recall that $\alpha$ is an atom of a measure $\nu$ if and only if
$F_\nu(\alpha)=0$ and the Julia-Carath\'{e}odory derivative
$F_\nu'(\alpha)\in[1,\infty)$, in which case
$\nu(\{\alpha\})=1/F_\nu'(\alpha)$. The atoms of $\mu^{\uplus q}$,
$q>0$, are characterized in the following proposition, which is a
direct consequence of Lemma \ref{basic} and the identity
$F_{\mu^{\uplus q}}'=qF_\mu'+1-q$, where $F_{\mu^{\uplus q}}'$ and
$F_\mu'$ are the Julia-Carath\'{e}odory derivatives.

\begin{prop} \label{qatom}
If $\mu\in\mathcal{M}$ $(\mu\neq\delta_a)$, $q>0$, and
$\alpha\in\mathbb{R}$ then $(1)$-$(3)$ are equivalent:
\begin{enumerate} [$\qquad(1)$]
\item {the point $\alpha$ is an atom of the measure $\mu^{\uplus
q}$;}
\item {$F_\mu(\alpha)=\alpha/q^*$ and the Julia-Carath\'{e}odory derivative
$F_\mu'(\alpha)\in(1,\infty)$;}
\item {$F_\mu(\alpha)=\alpha/q^*$ and $f_\mu(\alpha)\in(0,\infty)$.}
\end{enumerate} If $r=(1-\mu^{\uplus q}(\{\alpha\}))^{-1}>1$ then
$(\ref{NeF})$ holds for $z=\alpha$ and
\[F_\mu'(\alpha)=1+f_\mu(\alpha)=1+\frac{1}{q(r-1)}.\]
\end{prop}

Using the identity $\mu=(\mu^{\uplus q})^{\uplus 1/q}$, $q>0$, gives
the following corollary.

\begin{cor} If $q>0$, $r>1$, and $\alpha\in\mathbb{R}$ then the following statements are
equivalent:
\begin{enumerate} [$\qquad(1)$]
\item {$\alpha$ is an atom of $\mu$;}
\item {$F_\mu(\alpha)=0$ and $f_\mu(\alpha)\in(0,\infty)$;}
\item {$F_{\mu^{\uplus q}}(\alpha)=(1-q)\alpha$ and the
Julia-Carath\'{e}odory derivative $F_{\mu^{\uplus
q}}'(\alpha)\in(1,\infty)$;}
\item {$F_{\mu^{\uplus1/q}}(\alpha)=\alpha/q^*$ and
the Julia-Carath\'{e}odory derivative
$F_{\mu^{\uplus1/q}}'(\alpha)\in(1,\infty)$.}
\end{enumerate}
If $\mu(\{\alpha\})=1-r^{-1}$ then $f_\mu(\alpha)=1/r-1$,
$F_{\mu^{\uplus q}}'(\alpha)=1+q/(r-1)$, and
$F_{\mu^{\uplus1/q}}'(\alpha)=1+[q(r-1)]^{-1}$.
\end{cor}

Next, we characterize the points in $\mathbb{R}$ at which $F_\mu$ is
defined, takes real values, and has finite Julia-Carath\'{e}odory
derivatives.

\begin{prop} \label{preal} Let $p>1$ and let $x,\alpha$, and $\beta$ be real
numbers. If $px+(1-p)\beta=\alpha$ then $(1)$-$(4)$ are equivalent:
\begin{enumerate} [$\qquad(1)$]
\item {$F_\mu(x)=\beta$ and $0<f_\mu(x)<1/(p-1)$;}
\item {$F_\mu(x)=\beta$ and Julia-Carath\'{e}odory derivative
$F_\mu'(x)\in(1,p^*)$;}
\item {$H_p(x)=\alpha$ and the Julia-Carath\'{e}odory derivative
$H_p'(x)\in(0,1)$;}
\item {$\omega_p(\alpha)=x$ and the Julia-Carath\'{e}odory derivative
$\omega_p'(\alpha)\in(1,\infty)$;}
\item {$F_{\mu^{\boxplus p}}(\alpha)=\beta$ and the Julia-Carath\'{e}odory derivative $F_{\mu^{\boxplus
p}}'(\alpha)\in(1,\infty)$.}
\end{enumerate} If $(1)$-$(5)$ holds
then $(2.1)$ holds for $z=x$ and
\[F_\mu'(x)=1+f_\mu(x)=\frac{p-H_p'(x)}{p-1}=\frac{p\omega_p'(\alpha)-1}{(p-1)\omega_p'(\alpha)}=\frac{pF_{\mu^{\boxplus
p}}'(\alpha)}{1+(p-1)F_{\mu^{\boxplus p}}'(\alpha)}.\]
\end{prop}

\begin{pf} The equivalence of (1) and (2) was proved in Lemma \ref{basic}. The equivalence of (2) and
(3) and that of (4) and (5) follow from the identities of
Julia-Carath\'{e}odory derivatives $H_p'(x)=p+(1-p)F_\mu'(x)$ and
$F_{\mu^{\boxplus p}}'(\alpha)=(p\omega_p'(\alpha)-1)/(p-1)$,
respectively. The implication that (4) implies (3) holds by the fact
$\omega_p'(\alpha)=1/H_p'(x)$, where the Julia-Carath\'{e}odory
derivative $\omega_p'(\alpha)$ must be understood as $+\infty$ if
the Julia-Carath\'{e}odory derivative $H_p'(x)=0$. Conversely, if
(3) holds, i.e., (1) holds (because (1) and (3) was proved to be
equivalent) then $x\in\overline{\Omega_p}$, whence we have
$\omega_p(H_p(x))=x$ by Proposition \ref{prop2.1} (3), and the
statement (4) holds by the equality $\omega_p'(\alpha)=1/H_p'(x)$.
The last assertion follows from the preceding discussions.
\end{pf} \qed

We are now in a position to characterize the atoms of the measures
in $\mathbb{B}_{p,q}(\mathcal{M})$.

\begin{prop} \label{pqatom} Suppose $\alpha\in\mathbb{R}$, $p>1$, and $q>0$, and let $p',q'$ be the numbers
defined in Proposition $\ref{3.6}$. If $p^*q\neq1$ then the
following statements are equivalent:
\begin{enumerate} [$\qquad(1)$]
\item {the point $\alpha$ is an atom of $\mathbb{B}_{p,q}(\mu)$;}
\item {$F_{\mu^{\boxplus p}}(\alpha)=\alpha/q^*$ and the
Julia-Carath\'{e}odory derivative $F_{\mu^{\boxplus
p}}'(\alpha)\in(1,\infty)$;}
\item {$F_\mu(\alpha/p')=\alpha/q^*$ and the
Julia-Carath\'{e}odory derivative $F_\mu'(\alpha/p')\in(1,p^*)$;}
\item {$F_\mu(\alpha/p')=\alpha/q^*$ and
$0<f_\mu(\alpha/p')<1/(p-1)$.}
\end{enumerate} Moreover, if
$\mathbb{B}_{p,q}(\mu)(\{\alpha\})=1-r^{-1}$ for some $r>1$ then
\[F'_{\mu^{\boxplus p}}(\alpha)=1+\frac{1}{q(r-1)}\] and
\[F_\mu'(\alpha/p')=1+f_\mu(\alpha/p')=\frac{pq(r-1)+p}{pq(r-1)+p-1}=1+\frac{1}{q'(rp'-1)}.\] In
addition, if $p^*q<1$ then $\mathbb{B}_{p,q}(\mu)$ has at most one
atom. Particularly, the above assertions hold for $\mathbb{B}_t$,
$t\in(0,\infty)\backslash\{1\}$, as well.
\end{prop}

\begin{pf} The equivalence of (1) and (2) follows from Proposition
\ref{qatom}. Next, note that the hypothesis $p^*q\neq1$ shows that
$p'\neq\infty$. Then letting $x=\alpha/p'$ and $\beta=\alpha/q^*$
gives the equivalence of (2) and (3) by Proposition \ref{preal}. By
Lemma \ref{basic} we see that (3) and (4) are equivalent. By simple
computations, the rest desired equalities also follow from Lemma
\ref{basic}, Proposition \ref{qatom}, and \ref{preal}. That the
measure $\mathbb{B}_{p,q}(\mu)$, $p^*q<1$, has at most one atom is a
direct consequence of Theorem \ref{3.3} and \ref{3.10}.
\end{pf} \qed

Proposition \ref{pqatom} indicates that the Julia-Carath\'{e}odory
derivative $F_\mu'<p^*$ is one of the necessary conditions to
guarantee the existence of an atom of $\mathbb{B}_{p,q}(\mu)$,
$p^*q\neq1$. Indeed, consider the symmetric Bernoulli distribution
$\mu=\frac{1}{2}(\delta_{-1}+\delta_1)$ and the arcsine law of
distribution $\mathbb{B}_{1/2}(\mu)$ whose density is given by
\[d(\mathbb{B}_{1/2}(\mu))(x)=\frac{1}{\pi\sqrt{2-x^2}}\;dx,\;\;\;\;\;[-\sqrt{2},\sqrt{2}].\]
In this case
$(p=3/2,q=2/3,p'=2,q'=1/2,p^*=3,\;\mathrm{and}\;q^*=-2)$, if
$\alpha=\pm\sqrt{2}$ then it is easy to check that
$F_\mu(\alpha/2)=-\alpha/2$ and the Julia-Carath\'{e}odory
derivative $F_\mu'(\alpha/2)=3$. However, the points $\pm\sqrt{2}$
fail to be atoms of $\mathbb{B}_{1/2}(\mu)$. This example also
reveals an inaccuracy in the statement of [Proposition 5.1(2),
\ref{BN1}], which only requires the Julia-Carath\'{e}odory
derivative $F_\mu'\leq p^*$.

Now we are ready to state the main theorem in this section whose
proof is basically based on Proposition \ref{prop2.1} and Theorem
\ref{Hthm}.

\begin{thm} \label{main} Suppose that $\mu$ is a measure $(\mu\neq\delta_a)$
in $\mathcal{M}$, and that $p>1,q>0$ such that $p^*q\neq1$
$(q'=1+pq-p\neq0)$. Using the notations in Proposition
$\ref{prop2.1}$ and Theorem $\ref{Hthm}$, the following statements
hold.
\begin{enumerate} [$\qquad(1)$]
\item {The nonatomic part of the measure $\mathbb{B}_{p,q}(\mu)$ is
absolutely continuous.}
\item {The measure $(\mathbb{B}_{p,q}(\mu))^{\mathrm{ac}}$ is concentrated on the closure
of $\psi_p(V_p^+)$.}
\item {The density of $\mathbb{B}_{p,q}(\mu))^{\mathrm{ac}}$ on the set $\psi_p(V_p^+)$ is given by
\[\frac{d(\mathbb{B}_{p,q}(\mu))^{\mathrm{ac}}}{dx}
(\psi_p(x))=\frac{(p-1)pqf_p(x)}{\pi|pqx-q'\psi_p(x)+ipqf_p(x)|^2}.\]}
\item {The density of $(\mathbb{B}_{p,q}(\mu))^{\mathrm{ac}}$ is
analytic on the set $\psi_p(V_p^+)$.}
\item {Let $n(p,q)$ be the number of the components in the support
of $(\mathbb{B}_{p,q}(\mu))^{\mathrm{ac}}$. Then $n(p_1,q_1)\geq
n(p_2,q_2)$ whenever $p_1\leq p_2$ and $q_1,q_2>0$.}
\end{enumerate}
Particularly, the statements $(1)$-$(5)$ holds for
$\mathbb{B}_t(\mu)$, $t\in(0,\infty)\backslash\{1\}$.
\end{thm}

\begin{pf} Since the function $\psi_p$ defined in Theorem \ref{Hthm} is a
homeomorphism on $\mathbb{R}$ and $\omega_p$ extends continuously to
$\mathbb{R}$ by Proposition \ref{prop2.1}(3), it follows from
(\ref{general}) that
\begin{equation} \label{F}
F_{\mathbb{B}_{p,q}(\mu)}(\psi_p(x))=\frac{pqx-q'\psi_p(x)+ipqf_p(x)}{p-1},\;\;\;\;\;x\in\mathbb{R}.
\end{equation}
Since $F_{\mathbb{B}_{p,q}}(\mu)$ extends continuously to
$\mathbb{C}^+\cup\mathbb{\mathbb{R}}$, by the inversion formula
(\ref{inversion}) we obtain
\[\frac{d(\mathbb{B}_{p,q}(\mu))^{\mathrm{ac}}}{dx}
(\psi_p(x))=\frac{(p-1)pqf_p(x)}{\pi|pqx-q'\psi_p(x)+ipqf_p(x)|^2},\;\;\;\;\;x\in
V_t^+.\] Comparing the above formula with (\ref{density})shows that
the supports of $(\mu^{\boxplus p})^{\mathrm{ac}}$ and
$(\mathbb{B}_{p,q}(\mu))^{\mathrm{ac}}$ coincide for any $q>0$.
Observe that $\Im\omega_p(\psi_p(x))=f_p(x)>0$ for $x\in V_p^+$,
whence $\omega_p$ is analytic on $V_p^+$ by Proposition
\ref{prop2.1}(4). From the preceding discussion, we deduce that
statements (2)-(5) hold by Theorem \ref{Hthm}.

Next, let $p'=pq/q'$. We claim that if a point $\alpha\in\mathbb{R}$
such that $F_{\mathbb{B}_{p,q}(\mu)}(\alpha)=0$ and the
Julia-Carath\'{e}odory derivative
$F_{\mathbb{B}_{p,q}(\mu)}'(\alpha)=\infty$ or, equivalently,
$F_\mu(\alpha/p')=\alpha/q^*$ and the Julia-Carath\'{e}odory
derivative $F_\mu'(\alpha/p')=p^*$, then $\alpha$ belongs to the set
$\psi_p\left(\overline{V_p^+}\right)$, which is the closure of
$\psi_p(V_p^+)$. Note that we have $f_p(\alpha/p')=0$ by Proposition
\ref{prop2.1}(2) and Lemma \ref{basic}, and there does not exist an
open interval $I$ containing $\alpha/p'$ such that $f_p(x)=0$ for
all $x\in I$. Indeed, if such an interval $I$ exists then
$\rho(I)=0$ by [Corollary 3.6, \ref{Huang}]. This implies that the
second order derivative of $f_\mu$ on $I$ is positive, whence
$f_\mu$ is strictly convex on $I$. But $f_\mu(x)\leq(p-1)^{-1}$ for
all $x\in I$ and $f_\mu(\alpha/p')=(p-1)^{-1}$, a contradiction.
This particularly implies that the point
$\alpha/p'\in\overline{V_p^+}$, whence
\[\psi_p(\alpha/p')=H_p(\alpha/p')=\frac{p\alpha}{p'}+(1-p)F_\mu(\alpha/p')=\alpha\in\psi_p
\left(\overline{V_p^+}\right),\] and the claim follows. Moreover, we
see that the set
\[\{x\in\mathbb{R}:f_p(x/p')=0,\;\;\psi_p(x/p')=x\;\;\mathrm{and}\;\;F_\mu'(x/p')<p^*\}\] is the
collection of all atoms of $\mathbb{B}_{p,q}(\mu)$ by Proposition
\ref{pqatom}, and the set
\[\{\psi_p(x/p'):x\in\mathbb{R},\;f_p(x/p')=0,\;\;\mathrm{and}\;\;\psi_p(x/p')\neq x\}\]
has $\mathbb{B}_{p,q}(\mu)$-measure zero by the established result
(2). Then the preceding discussions and the established results (2)
and (3) show that $\mathbb{R}=\psi_p(\mathbb{R})$ consists of the
atoms of $\mathbb{B}_{p,q}(\mu)$ and the support of
$(\mathbb{B}_{p,q}(\mu))^{\mathrm{ac}}$, and therefore the assertion
(1) follows.
\end{pf} \qed

For the rest of the paper, we turn the attention to numbers $p>1$
and $q>0$ such that $p^*q=1$. The following proposition follows from
Lemma \ref{basic}, Proposition \ref{qatom}, and \ref{preal} and the
proof is left to the reader.

\begin{prop} \label{p*} If $p>1$ and
$\alpha\in\mathbb{R}$ then the following statements are equivalent:
\begin{enumerate} [$\qquad(1)$]
\item {the point $\alpha$ is an atom of the measure
$\mathbb{B}_{p,1/p^*}(\mu)$;}
\item {$F_{\mu^{\boxplus p}}(\alpha)=\alpha/(1-p)$ and the Julia-Carath\'{e}odory
derivative $F_{\mu^{\boxplus p}}'(\alpha)\in(1,\infty)$.}
\item {$F_\mu(0)=\alpha/(1-p)$ and the Julia-Carath\'{e}odory
derivative $F_\mu'(0)\in(1,p^*)$;}
\item {$F_\mu(0)=\alpha/(1-p)$ and $0<f_\mu(0)<1/(p-1)$.}
\end{enumerate}
If $\mathbb{B}_{p,1/p^*}(\mu)(\{\alpha\})=1-r^{-1}$ for some $r>1$
then
\[F_\mu'(0)=1+f_\mu(0)=1+\frac{1}{r(p-1)}\;\;\;
\mathrm{and}\;\;\;F_{\mu^{\boxplus
p}}'(\alpha)=1+\frac{p}{(p-1)(r-1)}.\] Particularly, the above
statements also holds for $\mathbb{B}_1$.
\end{prop}

If $\mu_0$ is the measure defined in Proposition \ref{Poisson} then
it is clear that $F_{\mu_0}(0)=0$ and $f_{\mu_0}(0)=1$. This yields
that the compound free Poisson distribution $p(\lambda,\nu)$ has an
atom at $0$ of mass $1-\lambda$ for $0<\lambda<1$ and no atom for
$\lambda\geq1$ by Proposition \ref{p*}.

The following theorem is a reformulation of Theorem \ref{3.10} since
$\Omega=\Omega_p$, $\psi=\psi_p$, and $f=f_p$. Therefore, its proof
is practically identical with that of Theorem \ref{3.10} or Theorem
\ref{main}, and is omitted.

\begin{thm} \label{main1} If $\mu\in\mathcal{M}$ and
$p>1$ then the following statements hold.
\begin{enumerate} [$\qquad(1)$]
\item {The nonatomic part of $\mathbb{B}_{p,1/p^*}(\mu)$ is absolutely
continuous.}
\item {The measure $(\mathbb{B}_{p,1/p^*}(\mu))^{\mathrm{ac}}$
is concentrated on the closure of $\psi_p(V_p^+)$.}
\item {The density of $(\mathbb{B}_{p,1/p^*}(\mu))^{\mathrm{ac}}$ on the set $\psi_p(V_p^+)$ is given by
\[\frac{d(\mathbb{B}_{p,1/p^*}(\mu))^{\mathrm{ac}}}{dx}(\psi_p(x))=\frac{f_p(x)}{\pi(x^2+f_p^2(x))}.\]}
\item {The density of $(\mathbb{B}_{p,1/p^*}(\mu))^{\mathrm{ac}}$ is
analytic on the set $\psi_p(V_p^+)$.}
\item {The number of the components in the support of $(\mathbb{B}_{p,1/p^*}(\mu))^{\mathrm{ac}}$
is a decreasing function of $p$.}
\end{enumerate}
Particularly, the
above statements also holds for $\mathbb{B}_1$.
\end{thm}

Since $\mu=\nu^{\boxplus p}$ for some $\nu\in\mathcal{M}$ and $p>1$
if $\mathrm{Ind}(\mu)>0$ by Proposition \ref{3.5}, we have the next
result by Theorem \ref{main} and \ref{main1}.

\begin{cor} If $\mu\in\mathcal{M}$ with $\mathrm{Ind}(\mu)>0$ then $(\mu^{\uplus
q})^{\mathrm{ac}}$ and $\mu^{\mathrm{ac}}$ contain the same number
of components in their supports for any $q>0$.
\end{cor}

It was shown in [\ref{Huang}] that there exists a measure
$\mu\in\mathcal{M}$ such that $\mu^{\boxplus p}$ contains infinitely
many components in the support for any $p>1$. Since
$(\mathbb{B}_{p,q}(\mu))^{\mathrm{ac}}$ and $(\mu^{\boxplus
p})^{\mathrm{ac}}$ have the number of components in their supports,
we have the following result.

\begin{prop} For any $t>0$, there exists a measure
$\mu_t\in\mathcal{M}$ such that $\mathrm{Ind}(\mu_t)=t$ and the
support of $\mu_t$ contains infinitely many components.
\end{prop}

\section*{Acknowledgements} The author wishes to thank his advisor,
Professor Hari Bercovici, for his generosity, and invaluable
discussion during the course of the investigation.

\end{document}